\documentclass{amsart}
\newtheorem{theorem}{Theorem}[section]
\newtheorem{proposition}[theorem]{Proposition}
\newtheorem{lemma}[theorem]{Lemma}

\newtheorem{corollary}[theorem]{Corollary}
 
\newtheorem{remark}[theorem]{Remark} 

\newcommand{\kerr}{\mbox{Ker} }
\newcommand{\cokerr}{\mbox{Coker} }
\newcommand{\so}{\mathfrak{s}_0}
\newcommand{\s}{\mathfrak{s}}
\renewcommand{\r}{\mathfrak{r}}
\renewcommand{\t}{\mathfrak{t}}
\newcommand{\tk}{\mathfrak{t}}

\newcommand{\imm}{\mbox{Im} }

\newcommand{\zee}{\mathbb{Z}}
\newcommand{\cue}{\mathbb{Q}}

\newcommand{\arr}{\mathbb{R}}
\newcommand{\cee}{\mathbb{C}}

\newcommand{\spinc}{\mbox{spin${}^c$ }}
\newcommand{\cfhat}{\widehat{CF}}
\newcommand{\hfhat}{\widehat{HF}}
\newcommand{\hfkhat}{\widehat{HFK}}
\newcommand{\aalpha}{{\mbox{\boldmath $\alpha$}}}
\newcommand{\bbeta}{\mbox{\boldmath $\beta$}}
\newcommand{\ggamma}{\mbox{\boldmath $\gamma$}}
\newcommand{\x}{{\mbox{\bf x}}}
\newcommand{\y}{{\mbox{\bf y}}}
\renewcommand{\hat}{\widehat}
\newcommand{\SL}{{\mathfrak{sl}}}
\newcommand{\VZ}{V_{\mathbb Z}}
\renewcommand{\P}{{\mathcal P}}
\newcommand{\B}{{\mathcal B}}
\newcommand{\T}{{\mathcal T}}

\newcommand{\tstar}{\tilde{\star}}
\newcommand{\Sym}{\mbox{Sym}}

\newcommand{\ts}{\textstyle}
\newcommand{\tP}{\tilde{P}}
\newcommand{\hmbar}{\overline{HM}}
\newcommand{\OS}{{Ozsv\'ath-Szab\'o }}
\newcommand{\OnS}{{Ozsv\'ath and Szab\'o }}

\usepackage{amsfonts}
\usepackage{graphicx}
\usepackage{amssymb}
\usepackage{amsmath}
\usepackage{amscd}

\begin{document}

\title{On the Heegaard Floer homology of a surface times a circle}

\author{Stanislav Jabuka}

\author{Thomas E. Mark}

\maketitle

\begin{abstract}We make a detailed study of the Heegaard Floer homology of the product of a closed surface $\Sigma_g$ of genus $g$ with $S^1$. We determine $HF^+(\Sigma_g\times S^1,\s;\cee)$ completely in the case $c_1(\s) = 0$, which for $g\geq 3$ was previously unknown. We show that in this case $HF^\infty$ is closely related to the cohomology of the total space of a certain circle bundle over the Jacobian torus of $\Sigma_g$, and furthermore that  $HF^+(\Sigma_g\times S^1,\s;\zee)$ contains nontrivial 2-torsion whenever $g\geq 3$ and $c_1(\s) = 0$. This is the first example known to the authors of torsion in $\zee$-coefficient Heegaard Floer homology. Our methods also give new information on the action of $H_1(\Sigma_g\times S^1)$ on $HF^+(\Sigma_g\times S^1,\s)$ when $c_1(\s)$ is nonzero.
\end{abstract}

\section{Introduction}

\subsection{Preliminaries and motivation} 
The goal of this work is to give a calculation of the Heegaard Floer homology groups of the product $\Sigma_g\times S^1$ of a closed oriented surface $\Sigma_g$ of arbitrary genus $g\geq 0$ with a circle. This problem has been addressed in part in other places \cite{OSknot}, \cite{OSplumb}, mainly as illustration of calculational techniques (some of which are used here). Our interest is more particular. A standard construction in smooth and symplectic 4-dimensional topology is the normal connected sum, or ``fiber sum'' construction where, in the simplest instance, neighborhoods of closed surfaces of the same genus and trivial normal bundle are removed from two 4-manifolds and their complements glued along their boundary. Indeed, together with a short list of other basic constructions, this operation gives rise to all known examples of exotic closed 4-manifolds. These examples are typically distinguished via gauge-theoretic invariants using either ad-hoc arguments or an appeal to general product theorems (or a combination). The 4-manifold invariants of Ozsv\'ath and Szab\'o \cite{OS3,OS4} are expected to give the same information as the Seiberg-Witten invariants but are in some cases easier to work with 
(results in this direction include the vanishing theorems for sums along $L$-spaces by Ozsv\'ath and Szab\'o  \cite{OS3} and the rational blowdown formulae by Roberts   \cite{roberts}). An approach to understanding the behavior of \OS invariants under fiber sum is to interpret the invariants as the result of a certain pairing of relative invariants of the complements of the surfaces in each 4-manifold, invariants taking values in the Floer homology of the boundary $\Sigma_g\times S^1$. This strategy is carried out in \cite{us3}; the results of the latter article are largely based  on the techniques and theorems of the present paper. % can be seen as a first step toward such a product formula.

%but to date no general results have appeared giving the behavior of these invariants under operations more %interesting than (ordinary) connected sums (however, see \cite{roberts}). 

A more ambitious goal is to study the \OS invariants of general fibered 3-manifolds. There is a sense in which the Floer homology of such a 3-manifold $Y$ is a module over the Heegaard Floer homology of the trivial fibered 3-manifold $\Sigma_g\times S^1$ (and in particular $HF(\Sigma_g\times S^1)$ should have the structure of a ring). Indeed, the pair-of-pants construction should lead to a product $HF(M_f)\otimes HF(M_g) \to HF(M_{fg})$ where $f:\Sigma_g\to \Sigma_g$ is an orientation-preserving diffeomorphism and $M_f$ denotes the mapping torus of $f$ (c.f. \cite{salamon}). Again, the basic information relevant to the study of this structure is an understanding of the Floer homology of $\Sigma_g\times S^1$. 

This paper deals with Floer homology with ``untwisted'' coefficients, though much of the technical work accomplished here is applicable to the more general situation. For this reason, the results here are not directly relevant to the problems mentioned above. Indeed, the behavior of \OS invariants under composition of cobordisms is functorial only modulo an indeterminacy (requiring a sum over \spinc structures in general), which is best handled by passing to Floer homology with twisted coefficient modules. However, even the untwisted Floer homology of $\Sigma_g\times S^1$ is surprisingly rich for such a simple manifold, and we hope that the results here will have independent interest.

\subsection{Statement of results}
We summarize the results of the paper. First, as mentioned above, certain parts of this problem have been understood by Ozsv\'ath and Szab\'o. Recall that if $\s$ is a \spinc structure on $\Sigma_g\times S^1$ such that the Chern class $c_1(\s)$ is not a multiple of the Poincar\'e dual of the class $[pt\times S^1]\in H_1(\Sigma_g\times S^1;\zee)$, then $HF^+(\Sigma_g\times S^1,\s;\zee)$ vanishes by the adjunction inequality \cite{OS2}. Furthermore, the adjunction inequality also implies that if $c_1(\s)$ is dual to $2k[S^1]$ then we obtain a similar vanishing unless $|k|\leq g-1$ (or $k=0$ when $g=0$). (See section \ref{backgroundsec} below for a summary of the basic constructions and notation of Heegaard Floer homology.)

When $c_1(\s)$ is dual to $2k[S^1]$ with $k\neq 0$, \OnS proved in \cite{OSknot} that $HF^+(\Sigma_g\times S^1,\s;\zee)$ is isomorphic as a (relatively, cyclically) graded $\zee[U]$ module to the ordinary cohomology of the $d$-fold symmetric power of $\Sigma_g$, where $d = g-1-|k|$. Following their notation, we set
\[
X(g,d) = H^*(\mbox{Sym}^d(\Sigma_g);\zee)[-g] \cong \bigoplus_{i=0}^d \Lambda^iH^1(\Sigma_g;\zee)\otimes \T_0/(U^{i-d-1})
\]
as graded $\zee[U]$-modules; here $[-g]$ indicates that the grading has been shifted down by $g$, while $\T_n$ denotes the graded ring $\zee[U,U^{-1}]/U\cdot\zee[U]$, graded so that the smallest degree of a nonzero homogeneous element is $n$. The factor $\Lambda^iH^1(\Sigma_g;\zee)$ is declared to lie in degree $i-g$, and as usual $U$ carries grading $-2$.

\OnS also calculated $HF^+(\Sigma_g\times S^1, \s_0)$, where $\s_0$ is the {\spinc} structure with $c_1(\s_0) = 0$, in the cases $g=0$ \cite{OS1}, $g=1$ \cite{OS3}, and $g=2$ \cite{OSplumb}. Thus, at least at the level of graded $\zee[U]$-modules, all Heegaard Floer homology of $\Sigma_g\times S^1$ was previously known except $HF^+(\Sigma_g\times S^1, \s_0)$ with $g\geq 3$. The calculation of these groups is the main goal of this paper. 

To state the result we introduce some notation. If $V$ is a symplectic vector space of dimension $2n$ then there is a natural operation $\Lambda: \Lambda^*V\to \Lambda^{*+2}V$ given by wedge product with the symplectic form (after using the symplectic form to identify $V^*\cong V$), and a natural adjoint map $L: \Lambda^*V\to \Lambda^{*-2}V$. The {\it primitive} subspace $P^* = \bigoplus_j P^j$ of $\Lambda^*V $ is defined to be the kernel of $L$, where $P^j = \ker(L)\cap \Lambda^jV$, and the {\it coprimitive} subspace $\tP^* = \bigoplus_j\tP^j$ is the kernel of $\Lambda$. We have $P^j = 0$ unless $j\leq n$, while $\tP^j = 0$ unless $j\geq n$, and in fact $\Lambda^{n-j}: \Lambda^jV\to \Lambda^{2n-j}V$ is an isomorphism mapping $P^j\to \tP^{2n-j}$. In the following, we apply these facts to the vector space $V = H^1(\Sigma_g;\cee)$ with symplectic stucture induced by cup product.

\begin{theorem}\label{mainthm} For all $g\geq 1$ there is an isomorphism of graded $\zee[U]$-modules
\begin{eqnarray*}
HF^+(\Sigma_g\times S^1, \s_0;\cee) &\cong& HF^+_{red}(\Sigma_g\times S^1,\s_0;\cee) \oplus\\ && \hspace*{.25in}\bigoplus_{j\geq 0} (P^j\otimes \T_{-g+j+1/2}) \oplus \bigoplus_{j\geq g} (\tP^j\otimes \T_{-g+j-1/2}),
\end{eqnarray*}
where $P^j$ (resp. $\tP^j$) is the space of primitive (resp coprimitive) forms of degree $j$ in $\Lambda^*H^1(\Sigma_g;\cee)$ and is considered to lie in degree 0.

Furthermore, there is an isomorphim of graded abelian groups
\[
HF^+_{red}(\Sigma_g\times S^1,\s_0;\zee) \cong X_0(g,g-3)[\ts\frac{5}{2}],
\]
where $X_0(g,d)$ is the group underlying the $\zee[U]$-module $X(g,d)$ and we take $X(g,d) = 0$ if $d<0$.
\end{theorem}

Of course $\cee$ may be replaced by $\cue$ (or any field of characteristic $0$) where it appears above; we give the statement in this way because complex coefficients arise naturally in the calculations, in particular the $P^j$ and $\tP^j$ are complex vector spaces.

Some remarks are in order. First, note that the reduced group $HF^+_{red}(\Sigma_g\times S^1,\s_0;\zee)$ vanishes for $g\leq 2$ (this was observed already by Ozsv\'ath and Szab\'o), and is nontrivial for $g\geq 3$. In fact, the reduced group is nontrivial in degrees from $-g + \frac{5}{2}$ to $g-\frac{7}{2}$.

Second, the identification of $HF^+_{red}(\Sigma_g\times S^1,\s_0;\zee)$ with $X(g,g-3)[\frac{5}{2}]$ on the level of graded groups does not respect the action of $U$. See Remark \ref{remarks} below for more details.

The Floer homology in large degrees can be understood as follows. Once $d\geq g-\frac{1}{2}$, the groups $HF^+_d(\Sigma_g\times S^1,\s_0;\cee)$ are all isomorphic; according to the theorem they are identified with the sum $\bigoplus_j P^{2j}\oplus\tP^{2j-1}$ or $\bigoplus_j P^{2j-1}\oplus\tP^{2j}$ depending on the parity of $d-\frac{1}{2}$. Identifying the primitives and coprimitives as above, we see that in each sufficiently large degree there is an identification $HF^+_d(\Sigma_g\times S^1,\s_0;\cee)\cong P^*$ (forgetting the grading on the right hand side). A more suggestive interpretation can be given as follows.

Consider the Jacobian torus $T^{2g} = H^1(\Sigma_g; \arr/\zee)$ of $\Sigma_g$. Then $T^{2g}$ carries a natural ``symplectic'' cohomology class $\omega\in H^2(T^{2g};\zee)\cong \Lambda^2(H^1(\Sigma_g;\zee))$ corresponding to the intersection pairing on $\Sigma$. Let $E_g$ be the total space of the $S^1$ bundle $E_g\to T^{2g}$ whose Euler class is $e(E_g) = \omega$; then it follows quickly from the Gysin sequence that the cohomology $H^*(E_g;\cee)$ can be identified with the (co)primitive cohomology of $T^{2g}$: 
\[
H^j(E_g;\cee) \cong \left\{\begin{array}{ll} P^j & j\leq g \\ \tP^{j-1} & j\geq g+1\end{array}\right.
\]
Here ``primitive'' is defined via multiplication and contraction by $\omega$ on $H^*(T^{2g};\cee)\cong \Lambda^*H^1(\Sigma_g;\cee)$. The previous description of the large-degree behavior of $HF^+(\Sigma_g\times S^1,\s_0;\cee)$ can be summarized in the statement (which follows from Theorem \ref{mainthm})
\begin{equation}\label{HFinfty}
HF^\infty(\Sigma_g\times S^1,\s_0;\cee)\cong H^*(E_g;\cee)\otimes \zee[U,U^{-1}],
\end{equation}
after adjusting the grading on $HF^\infty$ by a half-integer. Note that it is not quite true that the image of $HF^\infty$ in $HF^+$ is $H^*(E_g;\cee)\otimes \zee[U,U^{-1}]/U\cdot\zee[U]$, since the coprimitive pieces appear in gradings too small by two.

Lee and Packer \cite{leepacker} calculated the integer cohomology $H^*(E_g;\zee)$ using the Gysin sequence together with results from combinatorial matrix theory. According to their results, $H^*(E_g;\zee)$ contains considerable torsion: indeed, cyclic summands of any given order $\leq n$ appear in $H^*(E_g;\zee)$ for all $g\geq 2n-1$. The calculations in the current paper that lead to \eqref{HFinfty} bear a substantial similarity to those in \cite{leepacker} (c.f. Remark \ref{contrrem} below), and lead one to suspect that \eqref{HFinfty} holds with $\cee$ replaced by $\zee$---in particular that torsion elements of any given order appear in $HF^+(\Sigma_g\times S^1,\s_0;\zee)$ for large enough $g$. At the moment, however, we are satisfied with the following.

\begin{theorem} The Floer homology $HF^\infty_d(\Sigma_g\times S^1,\s_0;\zee_2)$ has dimension $2^{2g-1} + 2^{g-1}$ for each $d$. 
\end{theorem}

Since $\dim HF^\infty_d(\Sigma_g\times S^1,\s_0;\cee) = \dim P^* = {{2g+1}\choose{g}}$ and $HF^+_d\cong HF^\infty_d$ for all sufficiently large $d$, we conclude:

\begin{corollary} The Floer homology $HF^\infty(\Sigma_g\times S^1,\s_0;\zee)$ (and therefore also $HF^+(\Sigma_g\times S^1,\s_0;\zee)$) contains $2$-torsion for all $g\geq 3$. In fact, the identification \eqref{HFinfty} also holds with $\cee$ replaced by $\zee_2$.
\end{corollary}

These results are proved in section \ref{mod2sec}, where the mod-2 Floer homology $HF^+(\Sigma_g\times S^1,\s_0;\zee_2)$ is determined. It can also be shown by direct (though computer-assisted) calculation using the approach in this paper that $HF^\infty(\Sigma_g\times S^1,\s_0;\zee)$---and therefore also $HF^+(\Sigma_g\times S^1,\s_0;\zee)$---contains elements of order 3 for all $g\geq 5$ and elements of order 4 for all $g\geq 7$. We remark that the finite version of Heegaard Floer homology, $\hfhat (\Sigma_g\times S^1,\s_0;\zee)$, is free abelian for all $g$ (c.f. Theorem \ref{hfhatanswer}).

The results of this paper give new information also in the case of {\spinc} structures whose Chern class is nontrivial. Recall that for any closed \spinc 3-manifold $Y$ there is an action of $\Lambda^*(H_1(Y;\zee)/\mbox{\it torsion})$ on the Heegaard Floer homology groups (where $H_1(Y;\zee)$ acts via a map of relative degree $-1$), written $\xi\mapsto A_\gamma(\xi)$ for $\xi \in HF^+(Y)$ and $\gamma\in H_1(Y)$. In \cite{OSknot}, \OnS determined this action of $H_1$ on $HF^+(\Sigma_g\times S^1,\s_k) \cong X(g,d)$, where $\s_k$ is the \spinc structure with $c_1(\s_k)$ dual to $2k[S^1]$, under the assumption $3|k|>g-2$. To state our results, the {\it standard} action of an element $\gamma\in H_1(\Sigma_g\times S^1)$ on a homogeneous element $\xi = x \otimes U^j\in X(g,d)$ is written $\gamma.\xi$, and is defined as follows. First, $\gamma.\xi = 0$ unless $\gamma\in H_1(\Sigma_g)\subset H_1(\Sigma_g\times S^1)$. In the latter case, we set
\[
\gamma.\xi = \gamma.(x \otimes U^j) = \iota_\gamma(x)\otimes U^j + PD(\gamma)\wedge x\otimes U^{j+1},
\]
where $\iota_\gamma$ denotes interior multiplication by $\gamma$. According to \cite{OSknot}, the action of $H_1(\Sigma_g\times S^1)$ on $HF^+(\Sigma_g\times S^1,\s_k)$ is this standard action when $3|k|>g-2$. It follows from Theorem \ref{Znontorsresult} below that in general there are corrections to this action, which take the form 
\[
A_\gamma(\xi) = \gamma.\xi + \rho_1(\gamma, \xi) + \rho_2(\gamma,\xi) + \cdots
\]
Here the $\rho_\ell$ represent the homogeneous terms (with respect to the absolute grading on $X(g,d)$ given previously) in a correction that may appear when $3|k|\leq g-2$. In the following we implicitly use the absolute grading on $X(g,d)$ described above, which lifts the relative $\zee/2|k|\zee$-grading on $HF^+(\Sigma_g\times S^1, \s_k)$.

\begin{theorem} \label{h1actiontheorem}
Identify $HF^+(\Sigma_g\times S^1,\s_k)$ with $X(g,d)$ as above, and let $\xi\in X(g,d)$ be an element of degree $n$. Then for $\gamma\in H_1(\Sigma_g)$, the degree of $\rho_\ell(\gamma,\xi)$ is $n - 1 - 2\ell |k|$. Furthermore, 
\[
\rho_\ell (\gamma,\xi) = 0 \quad\mbox{ if $n < (2\ell -1)|k|$.}
\]
More specifically, if $\xi = x\otimes U^q\in X(g,d)$ with $x\in \Lambda^pH^1(\Sigma_g)$, $p = n+g+2q$, then for each $\ell$ with $0 < \ell \leq [\frac{n+|k|}{2|k|}]$, the correction $\rho_\ell (\gamma,\xi)$ lies in $\Lambda^{a+2\ell |k| + 1}H^1(\Sigma_g)\otimes U^{b + 2\ell |k| + 1}$ where $a = g-2-2|k| - n$ and $b = -|k|-1-n$. 
\end{theorem}

The corrections $\rho_\ell (\gamma,\xi)$ not forced to vanish by the above result are typically nonzero; a recipe for calculating them is outlined in section \ref{nontorssec}.

One can understand the $H_1$-action diagrammatically as follows, (c.f. figure \ref{pic1} below). \begin{figure}[b]
\centering
\includegraphics[width=10cm]{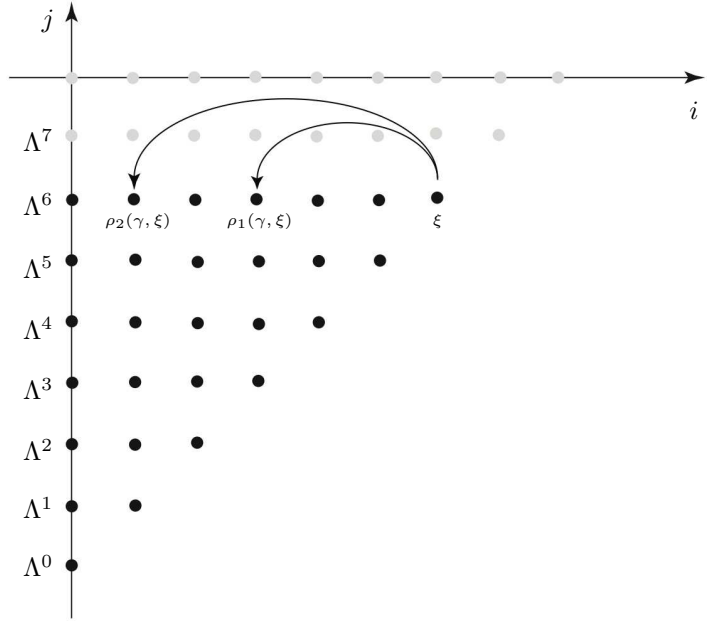}
\put(-15,215){$i$}
\put(-260,250){$j$}
\put(-112,175){\tiny $\xi$}
\put(-190,175){\tiny $\rho_1(\gamma, \xi)$}
\put(-236,175){\tiny $\rho_2(\gamma, \xi)$}
\put(-267,42){$\Lambda^0$}
\put(-267,64){$\Lambda^1$}
\put(-267,86){$\Lambda^2$}
\put(-267,109){$\Lambda^3$}
\put(-267,132){$\Lambda^4$}
\put(-267,155){$\Lambda^5$}
\put(-267,179){$\Lambda^6$}
\put(-267,203){$\Lambda^7$}
\caption{ A visualization of $X(g,d)$ and of the action $A_\gamma$ for the case of $g=8$ and $d=6$. 
The direct sum of the terms represented by the black dots yields $X(g,d)$. Each black dot along the $j$-axis represents an exterior power $\Lambda ^k =\Lambda ^k H^1(\Sigma _g ;\mathbb{Z})$ as indicated, the other black dots are translates of these along lines of slope 1. The amount of the translation is kept track of by a power of $U$, e.g. the dot at coordinates $(3,-4)$ represents $\Lambda ^1\otimes U^{-3}$. Taking, for example, $\xi \in \Lambda ^0 \otimes U^{-6}$, Theorem \ref{h1actiontheorem} shows that the expression $A_\gamma(\xi) - \gamma. \xi$ has two (potentially) nonvanishing ``correction terms'' $\rho_1(\gamma,\xi)$ and $\rho_2(\gamma,\xi)$, lying in $\Lambda^3\otimes U^{-3}$ and $\Lambda^{5}\otimes U^{-1}$ respectively. }  \label{pic1}
\end{figure}
We visualize $X(g,d)$ as a two dimensional array of groups, where the group at position $(i,j)$ is $\Lambda^p\otimes U^{-i}$ and $p=g-i+j$. Thus $X(g,d)$ lies in a triangle in the fourth quadrant of the $(i,j)$ plane, bounded by the lines $i = 0$, $j = d-g$, and $j = i-g$. The grading is given by $i+j$, so that $X(g,d)$ is symmetrical around degree $d-g$. The constraints given in the theorem indicate that nontrivial terms of $A_\gamma(\xi) - \gamma.\xi$ appear only in the top row of this triangle and in degrees $2\ell |k|$ less than the degree of $\gamma.\xi$. Observe that the largest degree of any element $\xi\in X(g,d)$ is $g-2|k|-2$, hence the degree of $\gamma.\xi$ is at most $g-2|k|-3$. The first correction to this appears in degree $(g-2|k|-3) - 2|k| = g-4|k|-3$, while the smallest degree of any element on the top row of $X(g,d)$ is $d-g = -|k|-1$. Thus if $g-4|k|-3 < -|k| - 1$ there can be no corrections, recovering the result of \OnS that the $H_1$-action is standard if $3|k|>g-2$.

At present, the methods employed in this article do not allow for a determination of the $\Lambda ^* H_1(\Sigma _g)$-module structure of $HF^+(\Sigma _g\times S^1,\s_k)$ in the case of $k=0$. As section \ref{mainsec} explains in greater detail, $HF^+(\Sigma _g \times S^1, \s_k)$ is computed by utilizing the surgery exact sequence 
of Ozsv\'ath and Szab\'o \cite{OS2}. In the case of $k\ne 0$, this approach identifies $HF^+(\Sigma _g \times S^1,\s_k)$ with the kernel of a map in this sequence, a kernel on which the $\Lambda ^*H_1(\Sigma _g)$-module structure is readily determined (and is as stated in Theorem \ref{h1actiontheorem}). However, in the case of $k=0$ 
$HF^+$ is a direct sum of a kernel and cokernel (as a group), which leads to a nontrivial extension problem for the module structure.

\subsection{Comparison to Seiberg-Witten-Floer homology} 

In \cite{km}, a work in preparation, Kronheimer and Mrowka
define three flavors of Seiberg-Witten Floer homology (or {\it monopole Floer homology}) for a spin$^c$ 3-manifold $(Y,\s)$, written
$$ \overline{HM}(Y,\s),  \quad \widehat{HM}(Y,\s), \quad \mbox{and}  \quad 
\stackrel{\rotatebox{180}{$\widehat{\phantom{.....}}$}}{HM}\hspace{-.5ex}(Y,\s).$$
A central conjecture in the theory states that these groups are 
isomorphic to $HF^\infty(Y,\s)$, $HF^-(Y,\s)$ and $HF^+(Y,\s)$ respectively. 

Kronheimer and Mrowka have obtained a general result on the structure of $\hmbar$, which can be described as follows.  Define a homomorphism $\beta_k:\Lambda ^kH^1(Y;\mathbb{Z}) \to \Lambda ^{k-3}H^1(Y;\mathbb{Z})$ by  
$$\beta_k(\alpha_1 \alpha_2...\alpha_k) = \sum _{i_1<i_2<i_3} (-1)^{i_1+i_2+i_3} \langle \alpha _{i_1} \cup \alpha _{i_2} \cup \alpha _{i_3} , [Y]\rangle  \,  \alpha _1 ... \hat{\alpha} _{i_1} ... \hat{\alpha} _{i_2} ...\hat{\alpha} _{i_3} ...{\alpha} _{k}.  $$
Then it is shown in section 35 of \cite{km} that there exists a decreasing filtration $\mathcal{F}_s$ on $\overline{HM}(Y,\s)$ whose 
associated graded groups are 
\begin{equation} \label{triplecup}
 \frac{\mathcal{F}_s \overline{HM}(Y,\s)}{\mathcal{F}_{s-1} \overline{HM}(Y,\s)} \cong 
\frac{\kerr (\beta_s)}{\imm (\beta_{s+3})} \otimes\zee[U,U^{-1}].
\end{equation}

Observe that the right hand side of \eqref{triplecup} is entirely determined by the triple-cup-product structure of $H^*(Y;\zee)$. In particular, cup products determine $\hmbar(Y,\s;\cue)$ as a rational vector space. Suitable analysis of the cup product structure on $\Sigma_g\times S^1$ yields an explicit description of $\hmbar(\Sigma_g\times S^1,\s_0;\cue)$ (proposition 35.3.4 of \cite{km}) which, though the approach is very different from the one taken here, is identical with our result for $HF^\infty(\Sigma_g\times S^1, \s_0)$ in \eqref{HFinfty}.

%Kronheimer and Mrowka compute $\overline{HM}$ (with rational coefficients) for the example of $Y = \Sigma _g \times S^1$ for all $g\ge 0$ and all spin$^c$-structures $\s$. Indeed, in the setting of SWF homology, 
%there is a purely topological description of $\overline{HM}(Y,\s)$ modulo torsion, in terms of triple cup products on $Y$, whenever $c_1(\s)$ is torsion. Namely, consider %the ring $\Lambda ^*H^1(Y;\mathbb{Z})$ and 
%the homomorphism $d_k:\Lambda ^kH^1(Y;\mathbb{Z}) \to \Lambda ^{k-3}H^1(Y;\mathbb{Z})$ given by  
%%
%$$d_k(\alpha_1 \alpha_2...\alpha_k) = \sum _{i_1<i_2<..<i_k} \langle \alpha _{i_1} \cup \alpha _{i_2} \cup \alpha _{i_3} , [Y]\rangle  \,  \alpha _1 ... \hat{\alpha} _{i_1} ... \hat{\alpha} _{i_2} ...\hat{\alpha} _{i_3} ...{\alpha} _{k}  $$
%%
%It was shown in \cite{km} that there exists a decreasing filtration $\mathcal{F}_s$ on $\overline{HM}(Y,\s)$ whose 
%associated graded groups are 
%%
%\begin{equation} \label{triplecup}
% \frac{\mathcal{F}_s \overline{HM}(Y,\s)}{\mathcal{F}_{s-1} \overline{HM}(Y,\s)} \cong 
%\frac{\kerr (d_s)}{\imm (d_{s+3})} 
%\end{equation}
%%
%The computations in \cite{km} of $\overline{HM}(\Sigma _g \times S^1,\s_0)$ are performed by evaluating the 
%right-hand side of this isomorphism. In particular, this approach is very different from the tools employed in the present article. Nevertheless, Kronheimer and Mrowka's results are in perfect agreement with our own, 
%thus providing new evidence towards the conjectured isomorphism between $\overline{HM}$ and $HF^\infty$. 

It is worth mentioning here that Ozsv\'ath and Szab\'o state a conjecture in \cite{OSplumb} regarding $HF^\infty(Y,\s)$ for a 3-manifold $Y$ with torsion \spinc structure $\s$, corresponding to the result of Kronheimer and Mrowka above. Specifically, there is a spectral sequence converging to $HF^\infty(Y;\s)$, whose $E_2$ and $E_3$ terms are isomorphic to $\Lambda^*H^1(Y;\zee)\otimes\zee[U,U^{-1}]$. \OnS conjecture that the $d_3$ differential in this sequence acts on $\Lambda^kH^1(Y;\zee)\otimes U^i$ by the homomorphism $\beta_k\otimes U^{-1}$, and that the higher differentials vanish---a conjecture that would imply an identification like \eqref{triplecup} for $HF^\infty$. The calculations in this paper can therefore be seen as evidence for their conjecture, and one might go so far as to attribute the appearance of torsion in the Heegaard Floer homology of $\Sigma_g\times S^1$ to the rich cup product structure on that manifold.

A slightly different definition of Seiberg-Witten Floer homology was considered by Mu\~noz and Wang in \cite{muwa}, and they give a calculation
for $\Sigma _g \times S^1$ in all non-torsion spin$^c$-structures. Their results agree with the calculations of $HF^+$ given in \cite{OSknot} and extended here.

\subsection{Organization}
The rest of the paper is organized as follows. In the next section we recall the relevant background material and notation for Heegaard Floer theory, all of which comes from \OnS \cite{OS1,OS2,OS3,OSknot}. That section outlines our calculational approach, which is based on knot Floer homology in much the same spirit as the calculations in \cite{OSknot}, \cite{us1,us2}. Section \ref{techsec} recalls the main tools from linear algebra that are necessary for the calculation, and sets up the preliminary results for section \ref{mainsec}. In the latter section the results listed above are proved.

{\bf Acknowledgements:} We are happy to thank Peter Kronheimer, Tom Mrowka, and Peter Ozsv\'ath for encouragement and useful comments. We are also indebted to the referee for a careful reading of this work and for many helpful remarks.

\section{Knot Floer homology and surgery}\label{backgroundsec}

We begin by recalling some of the relevant ideas and results of Heegaard Floer theory. The main point here is to outline an approach to the calculation of Heegaard Floer homology groups of a 3-manifold $Y_0$ obtained by zero-framed surgery on a nullhomologous knot $K$ in a closed, oriented 3-manifold $Y$. There is nothing essentially new here; we are making use of an idea introduced by Ozv\'ath and Szab\'o in \cite{OSknot} and by Rasmussen \cite{rasmussen}, and further exploited in \cite{us1} and \cite{us2}. In contrast to those situations, however, we are interested here particularly in the case of a \spinc structure on $Y_0$ whose first Chern class is a torsion element of $H^2(Y_0;\zee)$. 

If $Y$ is a closed oriented 3-manifold and $\s$ a \spinc structure on $Y$, we can find a pointed Heegaard diagram $(\Sigma, \aalpha,\bbeta,z)$ satisfying certain ``admissibility'' conditions with respect to $\s$. Here $\Sigma$ is a closed surface of genus $g$ embedded in $Y$ (the Heegaard surface), $\aalpha = \alpha_1,\ldots,\alpha_g$ and $\bbeta = \beta_1,\ldots,\beta_g$ are two $g$-tuples of pairwise disjoint essential simple closed curves on $\Sigma$, and $z$ is a point of $\Sigma$ away from the $\aalpha$ and $\bbeta$. In the symmetric power $Sym^g(\Sigma)$ we have two distinguished $g$-dimensional tori $T_\alpha = \alpha_1\times\cdots \times\alpha_g$ and $T_\beta = \beta_1\times\cdots\times\beta_g$ which intersect transversely, assuming the $\alpha_i$ intersect the $\beta_j$ transversely in $\Sigma$. There is a function $s_z: T_\alpha\cap T_\beta\to Spin^c(Y)$, depending on the basepoint $z$, that associates a \spinc structure to each intersection point $\x\in T_\alpha\cap T_\beta$. The Heegaard Floer chain complex $CF^\infty(Y;\s)$ is freely generated over $\zee$ by pairs $[\x,i]$, where $\x\in T_\alpha\cap T_\beta$ is an intersection point with $s_z(\x) = \s$ and $i$ is an arbitrary integer.

The boundary operator $\partial^\infty :CF^\infty\to CF^\infty$ ``counts'' certain holomorphic disks in $Sym^g(\Sigma)$ that connect pairs of intersection points, and has the property that any generator $[\y,j]$ appearing with nonzero coefficient in $\partial^\infty([\x,i])$ satisfies $j\leq i$. Therefore the subgroup $CF^-(Y,\s)$ generated by those pairs $[\x,i]$ with $i<0$ is a subcomplex, and the quotient group $CF^+(Y,\s)$ inherits the structure of a chain complex. There is a chain endomorphism $U: [\x,i]\mapsto [\x,i-1]$ of $CF^\infty$ that induces endomorphisms of $CF^-$ and $CF^+$; this map is of degree $-2$. The kernel of the action of $U$ on $CF^+$ is generated by pairs $[\x,0]$; the corresponding subcomplex of $CF^+$ is denoted $\cfhat$. Observe that as there are only finitely many intersection points $\x$, the complex $\cfhat$ is finitely generated and therefore the corresponding homology $\hfhat(Y,\s)$ is finitely generated as well. The main result of \cite{OS1} is that the homology groups of these various complexes, written
\[
HF^\infty(Y,\s),\quad HF^-(Y,\s),\quad HF^+(Y,\s),\quad \hfhat(Y,\s),
\]
are topological invariants of the pair $(Y,\s)$. 

Some of the properties of Heegaard Floer homology that we will use are as follows. Following standard practice we will write $HF^\circ$ to indicate that a given property holds for all of the above groups.
\begin{itemize}
\item There is a long exact sequence arising from the definition of $CF^+$ as the quotient of $CF^\infty$ by $CF^-$, namely
\[
\cdots \to HF^-(Y,\s) \to  HF^\infty(Y,\s)  \to  HF^+(Y,\s)  \to \cdots
\]
\item There are only finitely many \spinc structures $\s$ for which $HF^+(Y,\s)$ is nonzero; these ``basic classes'' are the same as those for which $\hfhat(Y,\s)$ is nontrivial in light of the following long exact sequence arising from the action of $U$:
\begin{equation}\label{Useq}
\cdots \to \hfhat(Y,\s) \to HF^+(Y,\s) \stackrel{U}{\to}  HF^+(Y,\s) \to \hfhat(Y,\s) \to \cdots 
\end{equation}
\item The homology groups $HF^\circ(Y,\s)$ are relatively cyclically graded, meaning that for a pair of generators $[\x,i]$, $[\y,j]$ there is a well-defined assignment of an element of $\zee/d(\s)\zee$ where $d(\s) = gcd_{A\in H_2(Y;\zee)}(\langle c_1(\s), A\rangle)$ is the divisibility of $c_1(\s)$.
\item In particular if $\s$ is a torsion \spinc structure---that is, $c_1(\s)$ is a torsion class---then $HF^\circ$ has a relative $\zee$ grading. In \cite{OS3} Ozv\'ath and Szab\'o ``lift'' this relative grading to an absolute, rational-valued grading.
\item If $W: Y_1\to Y_2$ is an oriented 4-dimensional cobordism between $Y_1$ and $Y_2$, equipped with a \spinc structure $\r$, then there are induced homomorphisms $F^\circ_{W,\r}: HF^\circ(Y_1,\s_1)\to HF^\circ(Y_2,\s_2)$ where $\s_i = \r|_{Y_i}$. If $\s_1$ and $\s_2$ are torsion \spinc structures, then $F^\circ_{W,\r}$ shifts absolute degree according to the following formula: for $\xi\in HF^\circ(Y_1,\s_1)$,
\begin{equation}\label{degformula}
\deg(F^\circ_{W,\r}(\xi)) - \deg(\xi) = \frac{1}{4}(c_1^2(\r) - 3\sigma(W) - 2\chi(W)),
\end{equation}
where $\sigma$ is the signature of the intersection form on $W$, $\chi(W)$ is the Euler characteristic, and $c_1^2(\r)$ makes use of the $\cue$-valued intersection pairing on elements of $H^2(W;\zee)$ whose boundary is torsion.
\item There is an action of the first homology $H_1(Y;\zee)$ on $HF^\circ(Y,\s)$ that decreases (relative) degree by one, written $\xi \mapsto \gamma.\xi$ for $\xi\in HF^\circ$ and $\gamma\in H_1(Y;\zee)$. In fact the torsion part of $H_1(Y;\zee)$ acts trivially, and the action extends to an action of $\Lambda^*(H_1(Y;\zee)/tors)$. The sequences and homomorphisms above are natural with respect to this action, and that of $U$. In particular, this means that if $\gamma_i\in H_1(Y_i;\zee)$ are homology classes on $Y_i$ for $i = 1,2$, and $\gamma_1$ is homologous to $\gamma_2$  in the cobordism $W:Y_1\to Y_2$, then $F^\circ_{W,\r}(\gamma_1.\xi) = \gamma_2. F^\circ_{W,\r}(\xi)$.
\end{itemize}

Now suppose that $K\subset Y$ is a nullhomologous knot. Oszv\'ath and Szab\'o explain in \cite{OSknot} how one can find a ``doubly-pointed'' Heegaard diagram $(\Sigma,\aalpha,\bbeta,w,z)$ in which $(\Sigma, \alpha_1,\ldots,\alpha_g, \beta_1,\ldots \beta_{g-1})$ is a Heegaard diagram for the complement of $K$ and $\beta_g$ is a meridian of $K$, while $w$ and $z$ are a pair of basepoints on either side of $\beta_g$. The additional basepoint provides a filtration of $CF^\infty(Y,\s)$, which is obtained as follows. By exchanging $\beta_g$ for a nullhomologous longitude $\gamma_g$ for $K$, we can obtain a Heegaard diagram $(\Sigma, \aalpha,\ggamma) = (\Sigma,\alpha_1,\ldots \alpha_g, \beta_1,\ldots,\beta_{g-1},\gamma_g)$ for the result $Y_0(K)$ of 0-framed surgery along $K$. An intersection point $\x\in T_\alpha\cap T_\beta$ gives rise in a certain canonical way to a pair of intersection points $\x_\pm\in T_\alpha\cap T_\gamma$ \cite{OSknot}. Fix a Seifert surface $S$ for $K$; then the filtration on $CF^\infty$ takes the form of an integer-valued function on generators where the ``filtration level'' of $[\x,i]$ is the integer $j$ given by
\[
j = i + \frac{1}{2}\langle c_1(s_z(\x_\pm)), [S_0]\rangle.
\]
Here $S_0$ denotes the closed surface in $Y_0(K)$ obtained by capping off the Seifert surface using the surgery disk, and the right-hand side is independent of the choice of $\x_+$ or $\x_-$. (In \cite{OSknot} the filtration is set up in slightly greater generality using relative \spinc structures.) The filtration is typically incorporated in the notation for generators of the Heegaard Floer complex by expressions of the form $[\x,i,j]$ (see figure \ref{pic2} below). 

\begin{figure}[b]
\centering
\includegraphics[width=12.5cm]{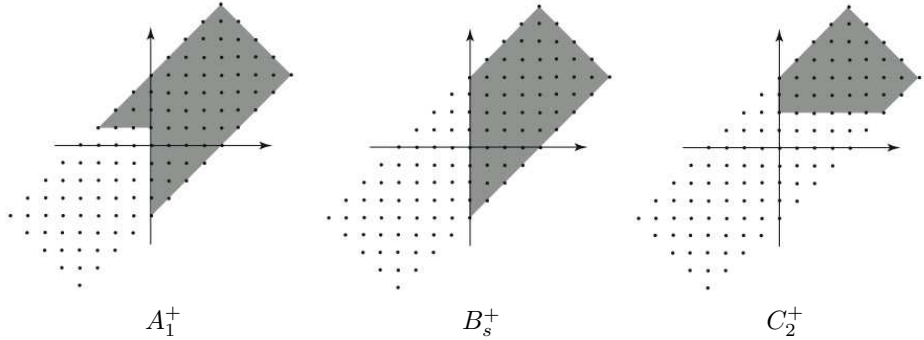}
\put(-300,-10){$A^+_1$}
\put(-180,-10){$B_s^+$}
\put(-65,-10){$C_2^+$}
\caption{A graphical representation of the complexes $A^+_s$, $B^+_s$ and $C^+_s$. In each of the three figures above, the dot at coordinates $(i,j)$ represents the free Abelian group generated by generators $[x,i,j]$. The complex $CFK^\infty$ is the direct sum over all such terms with $i,j\in \mathbb{Z}$. 
It follows from \cite{OSknot} that the nontrivial summands live in the region $i-g\le j\le i+ g$ where $g$ is the genus of the knot (thus the figures above correspond to $g=4$). The complexes $A_s^+$, $B_s^+$ and $C_s^+$ can be thought of as the direct sums of the terms in the shaded region (though they are quotients of $CFK^\infty$ by the direct sum of terms in the white regions). }  \label{pic2}
\end{figure}

Since the boundary operator is nonincreasing in both $i$ and $j$, we can specify subcomplexes of $CF^\infty$ by regions $R$ of the $(i,j)$ plane as follows: such a region corresponds to a subcomplex if whenever $(i,j)\in R$ we have $(i',j')\in R$ for all $i'\leq i$ and $j'\leq j$. The subcomplex corresponding to $R$ is generated by all elements $[\x,i,j]$ with $(i,j)\in R$, and is denoted by $C\{(i,j)\in R\}$. By abuse of notation, we use similar expressions to denote quotient complexes: thus $C\{i\geq 0\}$ is the quotient of $CF^\infty$ by the subcomplex $C\{i<0\}$ (and is isomorphic to $CF^+$ after forgetting the filtration). Of importance for our present purposes are the complexes
\[
A^+_s = C\{i\geq 0 \mbox{ or } j\geq s\} \qquad B^+ = C\{i\geq 0\} \qquad C^+_s = C\{i\geq 0 \mbox{ and } j\geq s\},
\]
for arbitrary $s\in \zee$. For a visual representation of these complexes, see figure \ref{pic2}. As just observed, $B^+$ is just $CF^+(Y,\s)$ with an additional filtration. The relevance of $A^+_s$ and $C^+_s$ is summarized in the following result, obtained by Ozsv\'ath and Szab\'o in \cite{OSknot}.

\begin{theorem}\label{OSsurgthm} For an integer $n\in \zee$, let $W_n$ denote the cobordism between $Y$ and $Y_n(K)$ obtained by adding a two-handle to $Y\times [0,1]$ along $K$ with framing $n$, and write $\overline{W}_n$ for $W_n$ with its natural orientation reversed, thought of as a cobordism from $Y_n(K)$ to $Y$. Write $\widehat{S}$ for the surface in $W_n$ obtained from the Seifert surface $S$ for $K$ union the core of the 2-handle. For a given $s\in \{0,\ldots,|n|-1\}$ and $\s\in Spin^c(Y)$ write $\t_s$ for the (unique) \spinc structure on $Y_n(K)$ cobordant to $\s$ via $W_n$ by a \spinc structure $\r_{s,k}$ satisfying $\langle c_1(\r_{s,k}), [\hat{S}]\rangle = (2k - 1) n - 2s$. Fix $s$. Then for all $n \gg 0$ there exists an isomorphism of chain complexes
\[
\Psi: CF^+(Y_n(K), \t_s) \to A^+_s.
\]
Furthermore, under the identification provided by $\Psi$, the chain map 
\[
f_{\overline{W}_n,\r_{s,0}}^+: CF^+(Y_n(K),\t_s) \to CF^+(Y,\s)
\] 
induced by $\overline{W}_n$ equipped with the \spinc structure $\r_{s,0}$ corresponds to the  natural projection $v^+: A^+_s\to B^+$.

Likewise, for all integers $-n\ll 0$, there is an isomorphism of chain complexes
\[
\Phi: C^+_s \to CF^+(Y_{-n}(K), \s_s)
\]
under which the homomorphism 
\[
f_{W_{-n},\r_{s,0}}^+: CF^+(Y,\s)\to CF^+(Y_{-n}(K),\t_s)
\]
corresponds to the natural projection $B^+\to C^+_s$.
\end{theorem}

In fact, it is observed in \cite{OSknot} that choosing $|n| \geq 2g(K)-1$ always suffices in the above theorem, where $g$ is the Seifert genus of $K$.

Now, since the Heegaard Floer homology of $Y$ can be calculated using either of the basepoints in the diagram $(\Sigma,\aalpha,\bbeta,w,z)$ associated with $K$, there exists a chain homotopy equivalence $J^+ : C\{i\geq 0\}\to C\{j\geq 0\}$. We will have occasion below to study this map in detail in a particular case; for now we merely recall the following, obtained in \cite{OSsurg}. Let $h^+: B^+ \to C_s^+$ denote the map given by the composition 
\[
h^+ = pr\circ U^{-s}\circ J^+
\]
where $pr: C\{j\geq s\}\to C\{i\geq 0 \mbox{ and } j\geq s\}$ is the obvious projection (the power of $U$ is included only for conceptual convenience, since $C\{j\geq s\}$ and $C\{j\geq 0\}$ are isomorphic as chain complexes). Then we have that under the identification provided by $\Psi$ above, the homomorphism
\[
f^+_{\overline{W}_{-n},\r_{s,1}}  :CF^+(Y,\s)\to CF^+(Y_{-n}(K),\t_s) 
\]
corresponds to $h^+$ (c.f. Theorem 2.3 of \cite{OSsurg}). 

With these preliminary notions established, we turn to the basic issue involved in our calculations: the determination of the Floer homology of the result of 0-surgery along a nullhomologous knot $K\subset Y$. 

The key tool here is the long exact sequence relating the Heegaard Floer homology of the manifolds $Y$, $Y_n$ and $Y_0$, for any integer $n\neq 0$. (For brevity we will write $Y_n$ in place of $Y_n(K)$ when the knot $K$ is understood from context.) Let $\s\in Spin^c(Y)$ be a torsion spin$^c$-structure on $Y$. Let $W_{-n}$ be the standard 2-handle cobordism from $Y$ 
to $Y_{-n}$ as before, and fix a Seifert surface $S\subset Y$ for $K$. Then we obtain an identification between $Spin^c(Y_0)$ and $Spin^c(Y)\times \mathbb{Z}$ where the projection to the
$\mathbb{Z}$ factor is given by $\mathfrak{u}\mapsto\langle c_1(\mathfrak{u}),[S_0] \rangle /2$.   
With this identification in place and with $\s \in Spin^c(Y)$ chosen, 
we let $[\s_s]$ denote the set
\[
[\s_s] = \{ (\s, s + \ell n ) \, | \, \ell\in \mathbb{Z} \}\subset Spin^c(Y_0),
\]
and write $HF^\circ(Y_0,[\s_s])$ for the direct sum of all Floer homologies in \spinc structures in $[\s_s]$. The surgery long exact sequence for any integer $n >0$ and any $s\in \{0,\ldots,n-1\}$ reads:
\begin{equation}\label{LESplus}
\cdots\to HF^+(Y_0,[\s_s])\stackrel{H}{\longrightarrow}HF^+(Y,\s)\stackrel{F}{\longrightarrow} HF^+(Y_{-n}, \t_s)\stackrel{G}{\longrightarrow}  HF^+(Y_0,[\s_s])\rightarrow\cdots
\end{equation}
where $\t_s\in Spin^c(Y_{-n})$ is as described in Theorem \ref{OSsurgthm}. There is a similar sequence with $HF^+$ replaced by $\hfhat$:
\begin{equation}\label{LEShat}
\cdots\to \hfhat(Y_0,[\s_s])\stackrel{H}{\longrightarrow} \hfhat(Y,\s)\stackrel{F}{\longrightarrow} \hfhat(Y_{-n}, \t_s)\stackrel{G}{\longrightarrow} \hfhat(Y_0,[\s_s])\rightarrow\cdots
\end{equation}
%
%In both sequences, the maps $F$, $G$ and $H$ are induced by counting holomorphic triangles (see \cite{OS3}).
%Moreover, the maps $G$ and $H$ lower the absolute grading by $1/2$. 

If $HF^\circ(Y,\s)$ is taken to be known, and the knot Floer homology of $K\subset Y$ is sufficiently well-understood that Theorem \ref{OSsurgthm} provides knowledge of $HF^\circ(Y_{-n},\t_s)$ (for large $n$), we see that a central issue in understanding $HF^\circ(Y_0,[\s_s])$ is the homomorphism $F$. By the construction of the long exact sequence $F$ is equal to the sum of maps induced by the cobordism $W_{-n}$ with the various \spinc structures  $\mathfrak{r}_{s,k} \in Spin^c(W_{-n})$ extending $\s$ and $\t_s$ (in the notation of Theorem \ref{OSsurgthm}).
Fix $s$ and let $F_k$ denote the component of $F$ corresponding to $\mathfrak{r}_{s,k}$. 
Since $\s$ and $\tk_s$ are torsion, the groups $\widehat{HF}(Y,\s)$ and $\widehat{HF}(Y_{-n},\tk_s)$
are equipped with an absolute $\mathbb{Q}$-grading. With respect to this grading, an easy calculation using \eqref{degformula} shows that the degree of $F_k$ is given by 
\begin{equation} \label{degreeshift}
\deg F_k = -nk^2 + (n-2s)k + \tau  \quad  \quad \mbox{ with } \quad \quad \tau = \frac{n-(2s-n)^2}{4n}
\end{equation}
Therefore the maximum degree shift of $F_k$ occurs when $k$ is the closest integer to $\frac{1}{2} - \frac{s}{n}$ (recall $s\in\{0,\ldots,n-1\}$), and in general the two components of highest degree in $F$ are $F_0$ and $F_1$. The components $F_k$ for $k\neq 0,1$ have degree at least $2n$ lower than that of $F_0$. (If $s = 0$, there are two components of $F$ realizing the maximal degree, namely those corresponding to $k= 0$ and $k = 1$.)
A similar calculation shows that $\deg H = -1/2$ and $\deg G = -\tau-1/2$ where $\tau$ is as in \eqref{degreeshift}. 

According to Theorem \ref{OSsurgthm} and the subsequent remarks, the maps induced by $F_0$ and $F_1$ are just $v^+$ and $h^+$ respectively, once $n$ has been chosen large enough that the theorem applies. Furthermore, if we consider only the case of $\hfhat$, we can choose $n$ so large that only $F_0$ and $F_1$ contribute to $F$ (indeed, we have seen that the degree shift associated to the other components of $F$ is roughly proportional to $n$: since $\hfhat$ is finitely generated, any map $\hfhat(Y,\s)\to\hfhat(Y_{-n},\t_s)$ that shifts degree sufficiently must be trivial). We have:

\begin{proposition}\label{hfhatcalc} Fix $n\gg 0$ and any $s\in\{0,\ldots,n-1\}$. Then under the identifications
\[
\hfhat(Y,\s) \cong H_*(C\{i= 0\}) \quad \mbox{ and } \quad \hfhat(Y_{-n},\t_s) \cong H_*(C\{\min\{i,j-s\} = 0\})
\]
the homomorphism $\widehat{F}: \hfhat(Y,\s)\to \hfhat(Y_{-n},\t_s) $ appearing in the surgery long exact sequence \eqref{LEShat} is identified with
\[
\widehat{F}  = \hat{v} + \hat{h}
\]
where $\hat{v}$ and $\hat{h}$ are the appropriate restrictions of $v^+$ and $h^+$.
\end{proposition}

The case of $HF^+$ is similar but slightly more involved: one makes use of a ``truncation'' argument analogous to those in \cite{OSsurg}. 

\begin{proposition}\label{hfpluscalc} Under the identifications
\[
HF^+(Y,\s) \cong H_*(C\{ i\geq0 \}) \quad  \mbox{ and } \quad  HF^+(Y_{-n},\tk_s) \cong H_*(C\{\min\{i,j-s\} \geq 0\})
\]
the homomorphism $F: HF^+(Y,\s) \to HF^+(Y_{-n},\t_s) $ can be identified with
\[
F^+ = v^+ + h^+.
\]
More precisely, there is an isomorphism of the Floer homology $HF^+(Y_0,\s_s)$ with the homology of the mapping cone of $v^++ h^+ : B^+\to C_s^+$.
\end{proposition}

\section{Preliminary Technicalities}\label{techsec}

\begin{figure}[b] 
\centering
\includegraphics*[width=12cm]{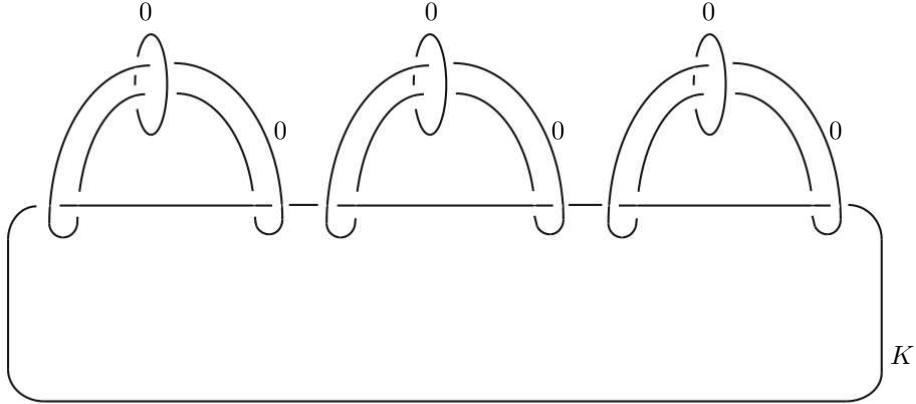}
\put(-2,30){$K$}
\put(-73,160){$0$}
\put(-178,160){$0$}
\put(-286,160){$0$}
\put(-25,115){$0$}
\put(-130,115){$0$}
\put(-235,115){$0$}
\caption{The knot $K$ in the case $g=3$.}  \label{pic3}
\end{figure}

We turn now to the case in which $Y_0(K) = \Sigma_g\times S^1$. Specifically, let $Y=\#^{2g}(S^1\times S^2)$ and let $K$ be the nullhomologous 
knot in $Y$ as indicated in figure \ref{pic3} (which  
depicts the case $g=3$), so that $Y_0 = \Sigma _g \times S^1$. Let $\s_0 \in Spin^c(Y)$ be the unique 
torsion spin$^c$-structure on $Y$.  
The knot Floer homology of $(Y,K)$ has been calculated in \cite{OSknot} (see also \cite{us1}) and equals
\begin{align} \label{knothomology}
\widehat{HFK}(Y,K,j,\s_0 ) & \cong \Lambda ^{g+j}H^1(\Sigma _g) \quad \mbox{ for } j=-g,...,g \quad \mbox{ and } \cr
HFK^\infty (Y,K,\s_0) & \cong \widehat{HFK}(Y,K,\s_0)\otimes \mathbb{Z}[U,U^{-1}]  
\end{align}
with the absolute grading on $\hfkhat$ agreeing with  the filtration level. The above is an isomorphism of 
$H_1(Y)\otimes \mathbb{Z}[U]$-modules where the action of $H_1(Y) \cong H_1(\Sigma _g)$ on the right-hand sides 
of  \eqref{knothomology} is given as 
\begin{align} \label{theaction}
\gamma . (x\otimes U^k)  = \iota _\gamma x \otimes U^k  + P.D. (\gamma) \wedge x \otimes U^{k+1}  
\quad \quad \gamma \in H_1(\Sigma_g), \, \, x\in \Lambda ^* H^1(\Sigma _g)\ 
\end{align}
where $\iota_\gamma x$ denotes the interior product of $x$ by $\gamma$, using the natural pairing between homology and cohomology.
 
To apply Proposition \ref{hfhatcalc} or \ref{hfpluscalc} to this situation, we must understand the projection $h^+ = pr\circ U^{s}\circ J^+$, and in particular the isomorphism $J^+: H_*(C\{i \geq 0\})\to H_*(C\{j \geq 0\})$. To express the result, it is convenient to have some additional background, which we introduce in the form of a digression.

\subsection{Digression on linear algebra}\label{linalgsec}
Suppose $V$ is a vector space of dimension $2g$ equipped with a symplectic form $\omega$. Throughout this section, we take the ground field to be $\mathbb C$. We can use $\omega$ to define a contraction $v\angle\cdot: \Lambda^kV\to \Lambda^{k-1}V$ for each $v\in V$ by the rule
\[
v\angle(w_1\wedge\cdots\wedge w_k) = \sum_\ell (-1)^{\ell-1} \omega(w_\ell, v) w_1\wedge\cdots\wedge \hat{w}_\ell\wedge\cdots\wedge w_k.
\]
This operation obeys the usual signed Leibniz rule, and extends to an action of $\Lambda^pV$ on $\Lambda^*V$: set 
\[
(v_1\wedge v_2\wedge\cdots \wedge v_p)\angle\alpha= v_1\angle(v_2\angle(\cdots (v_p\angle\alpha)\cdots)): \Lambda^kV\to \Lambda^{k-p}V.
\]
Furthermore, the symplectic form induces an identification $V\to V^*$ by $v\mapsto \omega(\cdot, v)$, and thence an isomorphism $\Lambda^*V\cong \Lambda^*V^*$. In what follows we will blur the distinction between $V$ and $V^*$, and in particular we can speak of contraction by the symplectic form $\omega\in \Lambda^2V^*$ itself.
\begin{lemma}\label{commlemma} Contraction by the symplectic form satisfies:
\[
\omega\angle(\omega\wedge\alpha) - \omega\wedge (\omega\angle\alpha) = (p-g)\alpha
\]
for $\alpha\in \Lambda^pV$.
\end{lemma}
\begin{proof} We can find a symplectic basis for $V$, i.e., a basis $\{e_i\}_{i=1}^{2g}$ with respect to which $\omega$ takes the standard form
\[
\omega(e_{2i-1},e_{2i}) = 1 = -\omega(e_{2i}, e_{2i-1}) \quad\mbox{and}\quad \omega(e_i, e_j) = 0 \mbox{ for other pairs of indices $i$, $j$}.
\]
Then we have (using juxtaposition to indicate wedge product):
\begin{eqnarray*}
\sum_i e_{2i} (e_{2i-1}\angle(e_{j_1}\cdots e_{j_p}))&=& \sum_{i,\ell}(-1)^{\ell-1} \omega(e_{j_\ell}, e_{2i-1},) e_{2i} e_{j_1}\cdots \widehat{e_{j_\ell}}\cdots e_{j_p}\\
&=& -n_{ev}\,e_{j_1}\cdots e_{j_p}
\end{eqnarray*}
where $n_{ev}$ denotes the number of even indices in $j_1,\ldots,j_p$. Similarly, 
\[
\sum_i e_{2i-1} (e_{2i}\angle (e_{j_1}\cdots e_{j_p}))= n_{odd}\, e_{j_1}\cdots e_{j_p},
\]
which gives the basic identity
\begin{equation}\label{pform}
\sum_i e_{2i-1}(e_{2i}\angle\alpha) - e_{2i} (e_{2i-1}\angle\alpha) = p\cdot\alpha \quad \mbox{for $\alpha$ a $p$-form}.
\end{equation}
Under our identification $V^*\cong V$, $\omega$ is given by $\sum_{i} e_{2i-1}\wedge e_{2i}$, and we easily find
\begin{equation}\label{omegacontr}
 \omega\angle\omega = \sum_{i,j = 1}^g  (e_{2i-1}e_{2i})\angle (e_{2j-1} e_{2j}) = -g.
\end{equation}
Now we calculate:
\begin{eqnarray*}
 \omega\angle(\omega\wedge\alpha) &=& \sum_i  (e_{2i-1} e_{2i})\angle (\omega\wedge\alpha)\\
&=& \sum_i  e_{2i-1}\angle( (e_{{2i}}\angle\omega)\wedge \alpha + \omega\wedge (e_{{2i}}\angle\alpha))\\
&=&  (\omega\angle\omega)\wedge\alpha +\\&& 
\left( \sum_i - (e_{{2i}}\angle\omega)\wedge (e_{{2i-1}}\angle\alpha) +  (e_{{2i-1}}\angle\omega)\wedge (e_{{2i}}\angle\alpha)\right) + \omega\wedge (\omega\angle\alpha)
\end{eqnarray*}
Since $ e_{k}\angle\omega = e_k$ for all $k$, the result follows from this, \eqref{pform}, and \eqref{omegacontr}.
\end{proof}

According to Lemma \ref{commlemma}, if we define operators $\Lambda$, $L$, and $H$ on $\Lambda^*V$ by
\[
\Lambda(\alpha) = \omega\wedge \alpha \qquad L(\alpha) =  -\omega\angle\alpha \qquad H(\alpha) = (p-g)\alpha
\]
for $\alpha\in \Lambda^pV$, then these three operators satisfy the relations
\[
[\Lambda, H] = -2\Lambda \qquad [L, H] = 2L \qquad [\Lambda, L] = H.
\]
In other words, the assignment $\left(\begin{array}{cc} 0 & 1 \\ 0 & 0 \end{array}\right) \mapsto \Lambda$, $\left(\begin{array}{cc} 0 & 0 \\ 1 & 0 \end{array}\right) \mapsto L$, and $\left(\begin{array}{cc} 1 & 0 \\ 0 & -1 \end{array}\right) \mapsto H$ endows $\Lambda^*V$ with the structure of an $\SL_2(\mathbb{C})$ representation. Appealing to the standard representation theory of $\SL_2(\mathbb{C})$, we obtain a decomposition of $\Lambda^{p}V$ (being the $(p-g)$-eigenspace of $H$) as a direct sum
\[
\Lambda^pV = P^p \oplus (\omega\wedge P^{p-2}) \oplus (\omega^2\wedge P^{p-4})\oplus\cdots,
\]
where $P^k = (\ker L)\cap \Lambda^kV$ is the ``primitive part'' of $\Lambda^kV$. Furthermore, $\Lambda^\ell = \omega^\ell\wedge\cdot : \Lambda^{g-\ell}V\to \Lambda^{g+\ell}V$ is an isomorphism for each $\ell$, respecting the above decomposition into primitive pieces.

To avoid combinatorial factors in later formulae, we define:
\[
\eta_k := \frac{\omega^k}{k!} = \sum_{j_1<\cdots < j_k} e_{2j_1-1}e_{2j_1}\cdots e_{2j_k - 1}e_{2j_k} \in \Lambda^{2k} V.
\]
Note that if $V$ is obtained from a symplectic integer lattice $\VZ$ by tensor product with $\mathbb C$ (for example, $V = H^1(\Sigma;\mathbb{C})$ and $\VZ = H^1(\Sigma; \mathbb{Z})$) then, supposing the basis $\{e_i\}$ to be induced from $\VZ$, the elements $\eta_k$ are defined over the integers. It is important to note, however, that the representation theory used above requires {\it complex} representations. In particular, we do not have a decomposition into sums of primitive factors over $\zee$.

The basic result that facilitates our computations is:

\begin{lemma}\label{swaplemma} For any $\xi\in \Lambda^p V$, we have 
\begin{equation}\label{contrswap}
\xi\angle \eta_k = \sum_{\ell \geq 0} (\eta_\ell\angle\xi)\wedge\eta_{k-p+\ell}.
\end{equation}
\end{lemma}

In particular, this shows that if $\beta\in P^q$ is primitive, then $\beta\angle\eta_\ell  = \beta\wedge \eta_{\ell - q}$.

\begin{proof}  By way of notation, let $z_i = e_{2i-1}\wedge e_{2i}\in \Lambda^2 V$. Easy calculations show that
\[
e_i\angle \eta_k = e_i\wedge\eta_{k-1}
\]
and
\[
z_i\angle\eta_k = -\eta_{k-1} + z_i\wedge \eta_{k-2}.
\]

Now, a basis for $\Lambda^pV$ is given by elements of the form $z_{i_1}\cdots z_{i_q}x_1\cdots x_s $, where juxtaposition indicates wedge product, $s + 2q =p$, and $x_1,\ldots,x_s\in V$ are members of the basis $\{e_1,\ldots,e_{2g}\}$, satisfying $\omega(x_i,x_j) = 0$ for all $i,j$.  We compute the contraction of $\eta_k$ by such an element:
\begin{eqnarray}
(z_{i_1}\cdots z_{i_q}x_1\cdots x_s)\angle\eta_k &=& x_1\cdots x_s ((z_{i_1}\cdots z_{i_q})\angle \eta_{k-s})\nonumber\\
&=& x_1\cdots x_s(( z_{i_1} \cdots z_{i_{q-1}})\angle (-\eta_{k-s-1} + z_{i_q}\eta_{k-s-2}))\nonumber\\
&=& x_1\cdots x_s ((z_{i_1}\cdots z_{i_{q-2}})\angle (\eta_{k-s-2} - (z_{i_q} + z_{i_{q-1}}) \eta_{k-s-3}\nonumber \\ && \hspace{5.65cm} + z_{i_q}z_{i_{q-1}} \eta_{k-s-4}))\nonumber\\
&=&  x_1\cdots x_s ((z_{i_1} \cdots z_{i_{q-3}})\angle ( - \eta_{k-s-3} \nonumber\\ &&\hspace{3.8cm} + (z_{i_q} + z_{i_{q-1}} + z_{i_{q-2}})\eta_{k-s-4} \nonumber\\ && \hspace{1.95cm} - (z_{i_q}z_{i_{q-2}} + z_{i_{q-1}}z_{i_{q-2}} + z_{i_q}z_{i_{q-1}})\eta_{k-s-5}\nonumber\\ && \hspace{4.7cm}+ (z_{i_q}z_{i_{q-1}}z_{i_{q-2}})\eta_{k-s-6}))\nonumber\\
&\vdots& \nonumber\\
&=& x_1\cdots x_s \sum_{r\geq 0} (-1)^{q-r} \sigma_r(z_{i_1},\ldots, z_{i_q}) \eta_{k-s-q-r},\label{leftside}
\end{eqnarray}
where $\sigma_r(z_{i_1},\ldots, z_{i_q})$ denotes the elementary symmetric polynomial of degree $r$ in the variables $z_{i_1},\ldots, z_{i_q}$.

On the other hand, it is easy to see that 
\[
z_j\angle ( x_1\cdots x_s z_{i_1}\cdots z_{i_q}) = \left\{\begin{array}{ll} -x_1\cdots x_s z_{i_1}\cdots \widehat{z_j}\cdots z_{i_q} & \mbox{if $j$ appears in $i_1,\ldots i_q$}\\
0 & \mbox{otherwise,}\end{array}\right.
\]
and since $\eta_n = \sigma_n(z_1,\ldots,z_g)$ we find
\[
\eta_n\angle(x_1\cdots x_s z_{i_1}\cdots z_{i_q}) = (-1)^n x_1\cdots x_s \sigma_{q-n}(z_{i_1},\ldots,z_{i_q}).
\]
Taking $\xi = x_1\cdots x_s z_{i_1}\cdots z_{i_q}$, the right hand side of \eqref{contrswap} becomes
\[
\sum_{\ell \geq 0} (\eta_\ell\angle\xi)\wedge\eta_{k-p+\ell} = 
x_1\cdots x_s\sum_{\ell= 0}^q (-1)^\ell \sigma_{q-\ell}(z_{i_1},\ldots,z_{i_q}) \eta_{k-p+\ell}
\]
Letting $r = q-\ell$ and recalling that $p = s+2q$, this is identical with \eqref{leftside}.
\end{proof}

Now consider the linear map $\tstar: \Lambda^*V\to \Lambda^{2g-*}V$ given by $\tstar: \alpha\mapsto \alpha\angle\eta_g$ for $\eta_g = \frac{1}{g!}\omega^g$ the volume form. This operator, analogous to the Hodge star associated to an inner product, will be called the {\it Hodge-Lefschetz star operator}. 

\begin{proposition}\label{starprop} For any $\beta\in P^q$ and any $n\geq 0$,
\[
\tstar (\beta\wedge\eta_n) = (-1)^{n}\beta\wedge\eta_{g-n -q},
\]
so that 
\[
\tstar^2(\beta\wedge\eta_n ) = (-1)^{g-q}\beta\wedge\eta_n. \qquad (\beta\in P^q).
\]

In particular, if $\beta\in P^{g-2k}$ then $\tstar (\beta\wedge\eta_k) = (-1)^{k} \beta\wedge \eta_k$. That is, under the decomposition
\[
\Lambda^gV = P^g \oplus (\omega\wedge P^{g-2}) \oplus (\omega^2\wedge P^{g-4})\oplus \cdots,
\]
the involution $\tstar: \alpha \mapsto  \alpha\angle\eta_g$ acts as $(-1)^{k}$ on the summand $\omega^k\wedge P^{g-2k}$.

\end{proposition}
Since the degree of $\beta\wedge \eta_n$ agrees modulo 2 with that of $\beta$ (that is, with $q$), we see that $\tstar^2$ acts as $(-1)^{g-p}$ on $\Lambda^{p}V$.

\begin{proof} Apply Lemma \ref{commlemma} to the product $\omega^m = \omega\wedge\omega^{m-1}$. We obtain
\begin{eqnarray*}
 \omega\angle\omega^m &=& \omega\wedge (\omega\angle\omega^{m-1}) - (g-(2m-2))\omega^{m-1}\\
&=& \omega^2\wedge (\omega\angle\omega^{m-2}) - [g-(2m-2) + g-(2m-4)]\omega^{m-1}\\
&\vdots&\\
&=& -[g-(2m-2) + g-(2m-4) + \cdots + g]\omega^{m-1}\\
&=& -m(g-m+1)\omega^{m-1}.
\end{eqnarray*}
Therefore,
\begin{eqnarray*}
 \omega^n\angle\omega^g &=& -g \omega^{n-1}\angle\omega^{g-1}\\
&=& (-g)(-(g-1)\cdot 2)) \omega^{n-2}\angle\omega^{g-2}\\
&\vdots&\\
&=& (-1)^n g(g-1)\cdots(g-n+1)\cdot n!\, \omega^{g-n}.
\end{eqnarray*}
Dividing by $n!g!$ we have
\[
 \eta_n\angle\eta_g = (-1)^n \eta_{g-n}.
\]
The proposition follows from this and the preceding lemma:
\[
(\beta\wedge \eta_n)\angle\eta_g = \beta\angle(\eta_n\angle\eta_g) = (-1)^n \beta\angle\eta_{g-n} = 
(-1)^{n} \beta\wedge \eta_{g-n-q}.
\]
\end{proof}

For later convenience we include the following result on the interaction between the Hodge-Lefschetz operator, wedge product, and contraction.

\begin{lemma}\label{stupidlemma} For any $v\in V$ and $\alpha\in \Lambda^*V$ we have the identities

1) $\tstar(v\wedge \alpha) =  v\angle\tstar \alpha$

2) $\tstar( v\angle\alpha) = -v\wedge \tstar\alpha$.
\end{lemma}

\begin{proof} The first equation is simply notation: $\tstar(v\wedge\alpha) =  (v\wedge\alpha)\angle\eta_g =  v\angle( \alpha\angle\eta_g) =  v\angle\tstar\alpha$. The second equation follows by applying the first to $\tstar\alpha$ and using $\tstar^2 = (-1)^{g-p}$ on $\Lambda^pV$.
\end{proof}

\begin{remark} The preceding is little more than a reformulation of the standard approach to the hard Lefschetz theorem (see, e.g., \cite{GH}) in the simple case of an exterior algebra of a complex vector space (or, more suggestively, the complex de Rham cohomology of a torus of real dimension $2g$). A slight difference is the omission here of a hermitian metric: the contraction is defined using the symplectic (K\"ahler) form itself. This makes a difference in the details---in particular, the behavior of the Hodge-Lefschetz operator here differs from that of the ordinary Hodge star---but it is easy to see that the primitive subspaces $P^k$ defined here are identical with those arising in the context of the hard Lefschetz theorem, where contraction typically refers to a hermitian inner product.
\end{remark}

This ends the digression on linear algebra.

\subsection{Calculations} We now return to the issue of understanding the isomorphism $J^+: H_*(C\{i\geq 0\})\to H_*(C\{j\geq 0\})$. Now, each of $H_*(C\{i\geq 0\})$ and $H_*(C\{j\geq 0\})$ is isomorphic as an abstract $\Lambda^*H_1(\Sigma_g)\otimes \zee[U]$ module to 
\[
M = \Lambda^*H^1(\Sigma)\otimes \T_0
\]
with the action of $H_1(\Sigma)$ given by \eqref{theaction}. It was observed by Ozsv\'ath and Szab\'o that the only module automorphisms of $M$ are multiplication by $\pm 1$ (section 5 of \cite{OSsurg}). 

\begin{figure}[b]
\centering
\includegraphics[width=14cm]{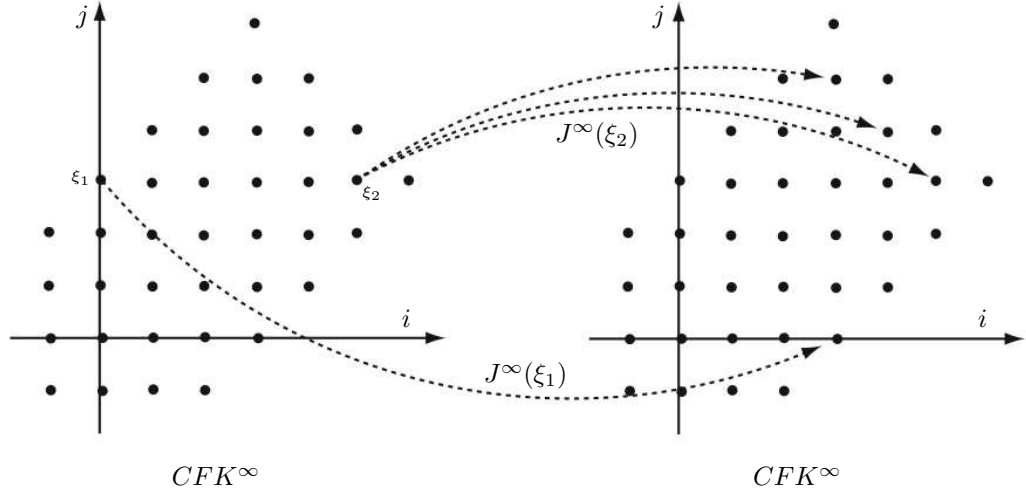}
\put(-25,55){$i$}
\put(-243,55){$i$}
\put(-148,170){$j$}
\put(-366,170){$j$}
\put(-330,-5){$CFK^\infty$}
\put(-110,-5){$CFK^\infty$}
\put(-212,34){$J^\infty(\xi _1)$}
\put(-185,125){$J^\infty (\xi_2)$}
\put(-368,110){\tiny $\xi_1$}
\put(-258,104){\tiny $\xi_2$}
\caption{A visualization of $J^\infty$. The two pictures above represent (portions of) two copies of $CFK^\infty$ where, as before, the dot at coordinates $(i,j)$ denotes the group $\Lambda ^{g+j-i}\otimes U^{-i}$ (with the above pictures showing the case of $g=3$). It is easily seen that the operator $\tstar$ maps terms at coordinates $(i,j)$ to those at coordinates $(j,i)$. On the other hand, the contraction operation $e^{2\omega U} \angle$ \lq\lq smears\rq\rq \, a term at coordinates $(i,j)$ downward along the line of slope -1. Thus, for example, $J^\infty(\xi_1)$ lives in $\Lambda ^0\otimes U^{-3}$ where 
$\xi _1 \in \Lambda ^6$ while $J^\infty (\xi_2)$ lives in $(\Lambda^5\otimes U^{-3})\oplus(\Lambda^3\otimes U^{-4})\oplus(\Lambda^1\otimes U^{-5})$ where $\xi_2\in \Lambda ^1\otimes U^{-5}$. 
}  \label{pic4}
\end{figure}

\begin{proposition}\label{Jprop} Under the usual identifications of $H_*(C\{i\geq 0\})$ and $H_*(C\{j\geq 0\})$ as submodules of $\Lambda^*H^1(\Sigma)\otimes \zee[U,U^{-1}]$, the action of $J^+$ is induced by the automorphism
\[
J^\infty(\xi\otimes U^i) = \varepsilon (-1)^p e^{2\omega U}\angle \tstar\xi \otimes U^j
\]
of $HFK^\infty(Y,K) = \Lambda^*H^1(\Sigma)\otimes\zee[U,U^{-1}]$, where $\omega$ denotes the intersection form on $H_1(\Sigma)$, thought of as an element of $\Lambda^2H^1(\Sigma)$. Here if $\xi\in\Lambda^pH^1(\Sigma)$ then we set $j = g+i-p$, and $\varepsilon = \pm 1$ is a sign depending only on $g$. Also, if $\alpha\in \Lambda^*H^1(\Sigma)$, we take $(\alpha\otimes U^n)\angle \xi$ to mean $(\alpha\angle \xi)\otimes U^{-n}$. 
\end{proposition}
For a pictorial description of $J^\infty$, consult figure \ref{pic4}. 

\begin{proof} Note that $e^{2\omega U} = \sum_{n\geq 0} 2^n\eta_n U^n$ is a $\zee[U,U^{-1}]$-linear combination of integer forms, having degree 0. It is straightforward to see that the homomorphism defined in the statement is an isomorphism of groups, and furthermore $J^\infty$ maps $C\{i<0\}$ to $C\{j<0\}$.  Thus $J^\infty$ induces an isomorphism $J^+: H_*(C\{i\geq 0\})\to H_*(C\{j\geq 0\})$. We will show below that $J^\infty$ respects the action of $H_1(\Sigma)$. Since, as mentioned above, the only 
$\Lambda ^*H_1(\Sigma _g)\otimes \mathbb{Z}[U]$-module isomorphisms of $\Lambda ^*H^1(\Sigma _g) \otimes 
\T_0$ are $\pm$id, the proposition follows.  

To this end, suppose $\gamma\in H_1(\Sigma)$, and observe that our sign conventions are arranged so that $\iota_\gamma(\cdot) = \gamma^*\angle\,\cdot\,$ (in the notation of the previous subsection), where $\gamma^*$ denotes the Poincar\'e dual of $\gamma$. Thus, taking $\varepsilon = +1$ for convenience,
\begin{eqnarray}
J^\infty(\gamma.\xi\otimes U^i) &=& J^\infty(\iota_\gamma\xi\otimes U^i + \gamma^*\wedge\xi\otimes U^{i+1})\nonumber\\
&=& (-1)^{p-1} e^{2\omega U}\angle(\tstar(\gamma^*\angle \xi)\otimes U^{j+1} + \tstar(\gamma^*\wedge\xi)\otimes U^j)\nonumber\\
&=& (-1)^{p-1} e^{2\omega U}\angle (-\gamma^*\wedge\tstar\xi\otimes U^{j+1} + \gamma^*\angle \tstar \xi\otimes U^j),\label{Jtemp}
\end{eqnarray}
using Lemma \ref{stupidlemma}. Now, it is not hard to check that $\omega\angle(\gamma^*\wedge\alpha) = \gamma^*\angle \alpha + \gamma^*\wedge(\omega\angle\alpha)$ for any $\alpha\in \Lambda^*H^1(\Sigma)$. It follows that 
\[
\omega^n\angle(\gamma^*\wedge \alpha) = n\gamma\angle(\omega^{n-1}\angle \alpha) + \gamma\wedge(\omega^n\angle \alpha)
\]
and thence
\[
e^{2\omega U}\angle(\gamma^*\wedge\alpha) = 2U^{-1}\gamma^*\angle(e^{2\omega U}\angle \alpha) + \gamma^*\wedge(e^{2\omega U}\angle \alpha).
\]
With this, \eqref{Jtemp} becomes
\begin{eqnarray*}
J^\infty(\gamma.\xi\otimes U^i) &=& (-1)^{p-1}(-e^{2\omega U}\angle(\gamma^*\wedge\tstar\xi)\otimes U^{j+1} + e^{2\omega U}\angle \gamma^*\angle\tstar\xi\otimes U^j) \\
&=& (-1)^{p-1}(-(2U^{-1} \gamma^*\angle e^{2\omega U}\angle\tstar\xi +\\&&\hspace*{.5in} \gamma^*\wedge(e^{2\omega U}\angle\tstar\xi))\otimes U^{j+1} + e^{2\omega U}\angle \gamma^*\angle\tstar\xi\otimes U^j )\\
&=& (-1)^{p}(\gamma^*\angle (e^{2\omega U}\angle\tstar\xi)\otimes U^{j} + \gamma^*\wedge(e^{2\omega U}\angle\tstar \xi)\otimes U^{j+1})\\
&=& \gamma.J^\infty(\xi\otimes U^i).
\end{eqnarray*}
\end{proof}

Note that while some of the results of the previous subsection hold only over $\cee$ (notably the primitive decomposition), it is easy to see that the results needed for the proof above and for the one to follow are true over $\zee$ as well.

In many situations, the sign $\varepsilon$ is unimportant (i.e., its value does not affect the result of the calculation of Floer homology groups). However, there are cases for which it is significant, and in section \ref{hatsubsec} below it will be shown that $\varepsilon = (-1)^{g-1}$. For this reason, we declare this to be the ``standard'' choice of sign, and formulate the following alternate description of the action of $J^\infty$ using it. 
\begin{lemma} \label{altJformulalemma}
For $\xi\in \Lambda^pH^1(\Sigma)$ and using the standard choice of sign, we have 
\begin{equation}\label{altJformula}
J^\infty(e^{-\omega U}\wedge\xi \otimes U^i ) = -e^{-\omega U} \wedge (e^{-\omega U}\angle\xi) \otimes U^i.
\end{equation}
\end{lemma} 

\begin{proof} By $U$-equivariance, it suffices to assume that $i = 0$. We calculate:
\begin{eqnarray*}
(-1)^{g-1}J^\infty(e^{-\omega U}\wedge\xi ) &=& (-1)^{g-1}\sum_{n= 0}^{g} J^\infty((-1)^n\eta_n\wedge\xi\otimes U^n)\\
 &=& \sum_{n= 0}^g (-1)^{p+j}e^{2\omega U}\angle\tstar(\xi\wedge\eta_n) \otimes U^{g-p-n}\\
&=& \sum_{n= 0}^g (-1)^{p+j}\sum_{q=0}^g 2^q\eta_q\angle \tstar(\xi\wedge\eta_n)\otimes U^{g-p-n-q}\\
&=& (-1)^p \sum_{j,q\geq 0} (-1)^j\xi\angle\tstar(2^q\eta_q\wedge\eta_n)\otimes U^{g-p-n-q}\\
&=& (-1)^p\xi\angle \,\sum_{m\geq 0}\sum_{q=0}^m (-1)^{m+q} \tstar(2^q\eta_q\wedge\eta_{m-q})\otimes U^{g-p-m}
\end{eqnarray*}
using Lemma \ref{stupidlemma}, and where $m = j+q$ in the last line. By expanding the identity $e^{\omega U} = e^{2\omega U}e^{-\omega U}$ and equating powers of $U$, we find 
\[
\sum_{q=0}^m (-1)^{m-q}2^q\eta_q\wedge\eta_{m-q} = \eta_m,
\]
so the preceding becomes
\begin{eqnarray*}
(-1)^{g-1}J^\infty(e^{-\omega U}\wedge\xi ) &=& (-1)^p\xi\angle\, \sum_{m\geq 0} \tstar\eta_m \otimes U^{g-p-m}\\
&=& (-1)^p\xi\angle\,\sum_{m\geq 0}(-1)^m\eta_{g-m}\otimes U^{g-p-m}\\
&=& \sum_{n\geq 0} (-1)^{g+p-n}\xi\angle\eta_n \otimes U^{n-p} \qquad \mbox{where $n = g-m$.}\\
\end{eqnarray*}
Applying Lemma \ref{swaplemma}, this becomes
\begin{eqnarray*}
\hspace*{1.5in} &=& \sum_{n\geq 0}(-1)^{g+p-n}\sum_{\ell\geq 0}(\eta_\ell\angle \xi)\wedge\eta_{n-p+\ell} \otimes U^{n-p}\\
&=& \sum_{r\geq p} (-1)^{g+p+r}\sum_{\ell=0}^r (-1)^\ell(\eta_\ell\angle\xi)\wedge\eta_{r-p}\otimes U^{r-\ell-p},
\end{eqnarray*}
where we let $r = n+\ell$ and note that $\eta_{r-p} = 0$ for $r<p$. Since $\eta_\ell\angle\xi = 0$ for $2\ell>p$, we may replace the upper limit on the $\ell$ summation by $g$ to obtain
\begin{eqnarray*}
\hspace*{1.5in} &=& \sum_{r\geq p} (-1)^{g+p+r}(e^{-\omega U}\angle \xi)\wedge\eta_{r-p}\otimes U^{r-p}\\
&=& (e^{-\omega U}\angle \xi)\wedge e^{-\omega U},
\end{eqnarray*}
which gives \eqref{altJformula}.
\end{proof}

\begin{corollary}\label{espacecor} For the standard sign choice, the $-1$ eigenspace of $J^\infty$ consists of all elements of the form $e^{-\omega U}\wedge \beta\otimes U^i$, where $\beta$ is primitive.
\end{corollary}

\begin{proof} From \eqref{altJformula}, $J^\infty(e^{-\omega U}\wedge\xi\wedge U^i) = -e^{-\omega U}\wedge\xi \otimes U^i$ if and only if
\[
e^{-\omega U}\wedge \left(\sum_{n\geq 1} (-1)^{n-1}\eta_n\angle \xi\right) = 0.
\]
Since multiplication with $e^{-\omega U}$ is an isomorphism of $\Lambda^*H^1(\Sigma)\otimes \zee[U,U^{-1}]$, this is equivalent to $\eta_n\angle\xi = 0$ for all $n>0$, i.e., $\xi$ is primitive.
\end{proof}

One can understand \eqref{altJformula} as saying that under the automorphism $\xi\mapsto e^{-\omega U}\wedge \xi$ of $\Lambda^*H^1(\Sigma)\otimes \zee[U,U^{-1}]$, the action of $J^\infty$ is identified with contraction by $-e^{-\omega U}$. This observation allows for relatively straightforward calculation of $HF^\infty(\Sigma\times S^1,\s_0)$ below.

Finally, $J^+$ and $\hat{J}$ are induced in the obvious way from $J^\infty$. In particular, since $U$ acts trivially on the ``hat'' complex we have that 
\[
\hat{J}: H_*(C\{i = 0\})\cong \Lambda^*H^1(\Sigma)\longrightarrow\Lambda^*H^1(\Sigma)\cong H_*(C\{j = 0\})
\]
is given by the ``top-order part'' of $J^\infty$, namely $\hat{J} = (-1)^{g+p-1}\tstar$ on $\Lambda^pH^1(\Sigma)$.

\section{Floer Homology for $\Sigma_g\times S^1$}\label{mainsec}

\subsection{Warmup: $\hfhat(\Sigma_g\times S^1, \s_0)$}\label{hatsubsec}
We focus here on the case of the torsion \spinc structure $\s_0$, since the other \spinc structures were dealt with in \cite{OSknot}.
Much of the answer follows easily from the form of the knot Floer homology of $K\subset \#^{2g} S^2\times S^1$, together with the surgery long exact sequence. Indeed, there is a spectral sequence associated to the filtration $[\x,i,j]\mapsto i+j$ of $CFK^{\infty}(Y,K)$ (or associated sub- or quotient complexes), and in the case at hand all differentials past the first are trivial. The $E_2$ term (and thus the $E_\infty$ term) is the appropriate sub- or quotient group of $HFK^\infty(Y,K)$, and in our case we see that for $n$ sufficiently large
\[
\widehat{HF}_{k+\tau}(Y_{-n}, \tk_0 ) \cong \left\{ 
\begin{array}{ll}
\Lambda ^{g+k} H^1(\Sigma _g ) \oplus [\Lambda ^{g-k} H^1(\Sigma _g )\otimes U^{-k}] \quad \quad & k\geq 1 \\
\Lambda ^g H^1(\Sigma _g ) & k=0 \cr
0 & \, k < 0
\end{array}
\right.
\]
Here $\tau$ is the degree shift \eqref{degreeshift} induced by the surgery cobordism.

From the surgery long-exact sequence \eqref{LEShat} we easily see 
$$ \widehat{HF}_{k + \frac{1}{2}}(Y_0,\so) \cong \widehat{HF}_{k} (Y,\s)\cong \Lambda^{g+k} H^1(\Sigma) \quad \quad \mbox{ for }
k \leq -2 $$
Since $Y_0$ admits an orientation reversing diffeomorphism we also get
\[
\widehat{HF}_{k - \frac{1}{2}}(Y_0,\so) \cong 
\widehat{HF}_{k} (Y,\s) \cong \Lambda^{g+k}H^1(\Sigma) \quad \quad \mbox{ for } k \ge 2 
\]
It remains to calculate $\widehat{HF}_{\frac{1}{2}}(Y_0, \so )$ (which is isomorphic to 
$\widehat{HF}_{- \frac{1}{2}}(Y_0, \so )$). For this we again use the surgery long exact sequence \eqref{LEShat}:
\begin{align} \nonumber
\cdots\to\widehat{HF}_{1}(Y,\s)  \stackrel{F_{(1)}}{\longrightarrow} \hfhat_{1+\tau}(Y_{-n}, \tk_0) \to \widehat{HF}_{\frac{1}{2}}(Y_0,\so) \rightarrow 
\widehat{HF}_{0}(Y,\s) \stackrel{F_{(0)}}{\longrightarrow} \widehat{HF}_{\tau}(Y_{-n},\tk_0) 
\rightarrow  \cdots 
\end{align}
Here $F_{(i)}$ denotes the restriction of $F$ to $\hfhat_{i}(Y,\s)$ (note that while $F$ is generally not homogeneous, for large $n$ we know that $F$ is given by its two top-degree parts and that when $s = 0$ these two parts have the same degree).

Now, according to Proposition \ref{hfhatcalc} and the observation at the end of the last section, we have (using the standard choice of sign $\varepsilon = (-1)^{g-1}$)
\[
F_{(i)}(\xi) = \left\{\begin{array}{ll} (\xi,\, (-1)^{i+1}\,\tstar\xi) \in \Lambda^{g+i}H^1(\Sigma)\oplus \Lambda^{g-i}H^1(\Sigma) & i\geq 1\\ \xi -\tstar\xi \in \Lambda^g H^1(\Sigma) & i = 0 \\ 0 & i< 0, \end{array}\right.
\]
where $\xi\in \hfhat_i(Y,\s) \cong \Lambda^{g+i}H^1(\Sigma)$. Therefore, in particular $F_{(1)}$ is injective with cokernel isomorphic to $\Lambda^{g+1}H^1(\Sigma)$. On the other hand, $\kerr (F_{(0)})$ is identified with the $+1$-eigenspace of $\tstar$ acting on $\Lambda^gH^1(\Sigma)$. In light of Proposition \ref{starprop}, this is:
\[
\kerr(F_{(0)})\cong {\mathcal P}_+ := P^g \oplus (\omega^2\wedge P^{g-4})\oplus (\omega^4\wedge P^{g-8})\oplus\cdots
\]
That is to say, $\kerr(F_{(0)})$ consists of ``half'' the primitive summands of $\Lambda^gH^1(\Sigma)$. In particular, since the rank of $P^k$ is ${{2g}\choose{k}} - {{2g}\choose{k-2}}$, it is not hard to see that the dimension of ${\mathcal P}_+$ as a vector space is $2^{g-1} + \frac{1}{2}{{2g}\choose{g}}$. (For example, apply de Moivre's formula $(1+i)^{2g} = 2^g ( \cos (g\pi/2) + i \sin (g\pi/2))$ and expand the left-hand side using the binomial theorem. We also give a direct calculation of this dimension after Theorem \ref{hfhatanswer}.) The results above can be summarized as follows.

\begin{proposition}\label{hfhatCanswer} For any half-integer $i\in \frac{1}{2} + \zee$ there is an isomorphism of vector spaces
\[
\hfhat_i(\Sigma\times S^1,\s_0; \cee) \cong \left\{\begin{array}{ll} \Lambda^{g- |i|-1/2}H^1(\Sigma; \cee) &  3/2 \leq |i| \leq g-1/2\\
\Lambda^{g-1}H^1(\Sigma; \cee)\oplus {\mathcal P}_+ & |i| = 1/2 \\ 0 & \mbox{otherwise} \end{array}\right.
\]
\end{proposition}

In fact we can be more precise: the Floer homology with $\zee$ coefficients is not difficult to obtain in this case. As a byproduct, we verify the sign choice $\varepsilon = (-1)^{g-1}$ used above. 

First, we claim that $\hfhat(\Sigma_g\times S^1,\s_0;\zee)$ is necessarily free abelian. To see this, observe that the surgery sequence above implies $\hfhat_i(\Sigma\times S^1,\s_0;\zee)\cong \Lambda^{g-|i|-1/2}H^1(\Sigma;\zee)$ when $|i|\geq \frac{3}{2}$, because for $i<-\frac{3}{2}$ we have an isomorphism $\hfhat_i(\Sigma\times S^1,\s_0;\zee)\to \hfhat_{i-1/2}(Y;\zee)$ and in general we know $\hfhat_i(\Sigma\times S^1,\s_0;\zee) \cong \hfhat_{-i}(\Sigma\times S^1,\s_0;\zee)$. To take care of the case $|i| = \frac{1}{2}$, observe that the surgery sequence gives us short exact sequences (with coefficients in $\zee$)
\begin{align*}
 0\to \cokerr(\Lambda^{g+1}\to \Lambda^{g+1}\oplus \Lambda^{g-1})\to \hfhat_{1/2}(\Sigma&\times S^1)\to \kerr(\Lambda^g\to \Lambda^g)\to 0 \\
 0\to \cokerr(\Lambda^g\to \Lambda^g)\to \hfhat_{-1/2}(\Sigma \times S^1&)\to \kerr(\Lambda^{g-1}\to 0)\to 0
\end{align*}
where the maps connecting exterior powers are given by $F$ in the appropriate degree. Certainly $\kerr(\Lambda^g\to \Lambda^g)$ is free; we also observe that the homomorphism $F_{(1)}: \Lambda^{g+1}\to \Lambda^{g+1}\oplus\Lambda^{g-1}$ given by $F_{(1)}(\xi) = (\xi, (-1)^{i+1}\varepsilon\tstar\xi)$ has free cokernel isomorphic to $\Lambda^{g+1}$ regardless of the value of $\varepsilon$. Therefore, $\hfhat_{1/2}(\Sigma\times S^1)$ is free abelian, and by duality we infer $\hfhat_{-1/2}(\Sigma\times S^1)$ is likewise free. This verifies the group structure of the Floer homology, and from this we can deduce that $\varepsilon = (-1)^{g-1}$. Indeed, from the exact sequence above it must be that $\cokerr(F_{(0)} = 1+(-1)^g\varepsilon\tstar: \Lambda^g\to \Lambda^g)$ is free; the latter can be calculated directly as follows.

Recall that a basis for $\Lambda^pH^1(\Sigma;\zee)$ is given by elements of the form $x_{i_1}\cdots x_{i_r}z_{j_1}\cdots z_{j_s}$ where
\begin{enumerate}
\item  $x_{i_n}$ is either $e_{2i_n-1}$ or $e_{2i_n}$ where $\{e_{2i-1},e_{2i}\}_{i=1}^g$ is a symplectic basis for $H^1(\Sigma;\zee)$, 
\item  $z_j = e_{2j-1}\wedge e_{2j}$, and 
\item the index sets $I = \{i_1\,\ldots, i_r\}$ and $J = \{j_1,\ldots, j_s\}$ are disjoint and satisfy $r + 2s = p$.
\end{enumerate}

Given indexing sets $I=\{i_1,...,i_r\}$ and $J=\{ j_1,..,j_s\}$ as above, let $K=\{k_1,...,k_t\}$ be the set of indices disjoint from $I\cup J$ and such that $I\cup J\cup K = \{1,...,g\}$. For the next calculation, observe that 
$\eta_g = z_{i_1}...z_{i_r}z_{j_1}...z_{j_s}z_{k_1}...z_{k_t}$. With this in mind we compute
$\tstar(x_{i_1}\cdots x_{i_r}z_{j_1}\cdots z_{j_s})$: 
%It is a simple direct calculation that 
\begin{align} \nonumber
\tstar(x_{i_1}\cdots x_{i_r}z_{j_1}\cdots z_{j_s}) & = (x_{i_1}\cdots x_{i_r}z_{j_1}\cdots z_{j_s}) \angle (z_{i_1}...z_{i_r}z_{j_1}...z_{j_s}z_{k_1}...z_{k_t})  \cr
& = (-1)^s (x_{i_1}\cdots x_{i_r}) \angle (z_{i_1}...z_{i_r}z_{k_1}...z_{k_t})  \cr
& = (-1)^s x_{i_1}\cdots x_{i_r} z_{k_1}\cdots z_{k_t}
\end{align}
In particular, elements of the form $x_{i_1}\cdots x_{i_g}\in \Lambda^g$ are fixed by $\tstar$. Therefore
\[
F(x_{i_1}\cdots x_{i_g}) = x_{i_1}\cdots x_{i_g} + (-1)^g\varepsilon\tstar(x_{i_1}\cdots x_{i_g}) = \left\{\begin{array}{ll} 0 &\mbox{if $\varepsilon = (-1)^{g-1}$} \\ 2x_{i_1}\cdots x_{i_g} & \mbox{if $\varepsilon = (-1)^g$}\end{array}\right.
\]
The requirement that $\cokerr(F)$ be free dictates that $\varepsilon = (-1)^{g-1}$.

\begin{theorem}\label{hfhatanswer} For any $g\geq 1$, there is an isomorphism of graded abelian groups
\[
\hfhat_i(\Sigma_g\times S^1,\s_0;\zee) \cong\left\{\begin{array}{ll} \Lambda^{g-|i|-1/2}H^1(\Sigma;\zee) & 3/2\leq |i|\leq g-1/2\\ \Lambda^{g-1}H^1(\Sigma;\zee)\oplus\P_+ & |i| = 1/2\\ 0 & \mbox{otherwise},\end{array}\right.
\]
where $\P_+$ is a free abelian group of rank $2^{g-1} + \frac{1}{2}{{2g}\choose{g}}$.\hfill$\Box$
\end{theorem}

From the above, we can identify $\P_+$ with the subgroup of $\Lambda^gH^1(\Sigma;\zee)$ spanned by self-dual elements:
\[
\P_+ = \mbox{span}\{ x_{i_1}\cdots x_{i_r}(z_{j_1}\cdots z_{j_{n}} + (-1)^nz_{j_{n+1}}\cdots z_{j_{2n}}) \}
\]
where $ \{i_1,\ldots,i_r,j_1,\ldots,j_{2n}\} = \{1,\ldots,g\}$ and $r + 2n = g$. To verify the rank of this space, observe that when $n = 0$ there are two choices for each $x_i$ giving rise to $2^g$ such basis elements. Since the total number of elements in the basis is ${2g}\choose{g}$ there are ${{2g}\choose{g}} - 2^g$ remaining, which are paired to form self-dual elements as above. Hence the rank of $\P_+$ is $2^g + \frac{1}{2}({{2g}\choose{g}} - 2^g) = 2^{g-1} + \frac{1}{2}{{2g}\choose{g}}$.

\subsection{The case of $HF^+(Y_0,\s_k)$}\label{nontorssec} Though the group structure of $HF^+(Y_0,\s_k)$ for $k \neq 0$ was determined by Ozsv\'ath and Szab\'o in \cite{OSknot}, we revisit the calculation here. Indeed, the results above detailing the form of $J^+$ allow for a description of the $\Lambda^*H_1(Y_0)\otimes \zee[U]$-module structure of the Floer homology groups, which was obtained in \cite{OSknot} only for a restricted range of $k$. Throughout this subsection, all homology (Floer and ordinary) is taken with coefficients in $\cee$.

We work with the negative surgery sequence as before, and make use of Proposition \ref{hfpluscalc}. We proceed by calculating the kernel and cokernel of $F: HF^+(Y,\s)\to HF^+(Y_{-n}, \tk_k)$. We can write
\[
F = pr \circ (1 + U^{-k} J^+)
\]
where $pr$ is the projection to $C^+_k = H_*(C\{i\geq 0 \mbox{ and } j\geq k\}) = HF^+(Y_{-n},\tk_k)$. For convenience, we suppose here that $k < 0$; the case $k> 0$ follows by the conjugation invariance of Heegaard Floer homology and the case $k = 0$ will require special attention.

\begin{lemma} The restriction of $F$ to $H_*(C\{i\geq 0\mbox{ and } j\geq k\})$ is surjective. In particular, $HF^+(Y_0, \s_k)$ is isomorphic to $\kerr(F)$.
\end{lemma}

\begin{proof} If $\xi\in H_*(C\{i\geq 0\mbox{ and } j\geq k\})$, then since $k<0$, $F(\xi) = \xi + pr(U^{-k}J^+(\xi))$ has highest-degree part equal to $\xi$. The lemma follows immediately.
\end{proof}

Following \cite{OSknot}, let us write $X(g,d)$ for the graded vector space
\begin{equation}\label{xgddef}
X(g,d) = \bigoplus_{i = 0}^d \Lambda^i H^1(\Sigma_g)\otimes \zee[U^{-1}]/U^{i-d-1},
\end{equation}
graded so that $U$ has degree $-2$ and $\Lambda^iH^1(\Sigma)$ is supported in degree $i-g$. Note that $X(g,d)$ can be identified as a vector space with the homology of $\Sym^d(\Sigma)$ (after a grading shift). There is also a natural embedding 
\[
\iota: X(g,d)\to HF^+(Y) =  \Lambda^*H^1(\Sigma)\otimes \zee[U,U^{-1}]/U\cdot\zee[U]
\]
whose image is $ H_*(C\{i\geq 0 \mbox{ and } j < k\})$ where $d = g-1 + k = g-1 - |k|$ (since $k< 0$). The latter identification induces an action of $H_1(\Sigma)$ on $X(g,d)$ (namely \eqref{theaction}), which we extend to an action of $H_1(\Sigma\times S^1)$ by requiring the extra circle to act trivially. We will refer to this action as the ``standard'' one.

\begin{figure}[bt]
\centering
\includegraphics[width=10cm]{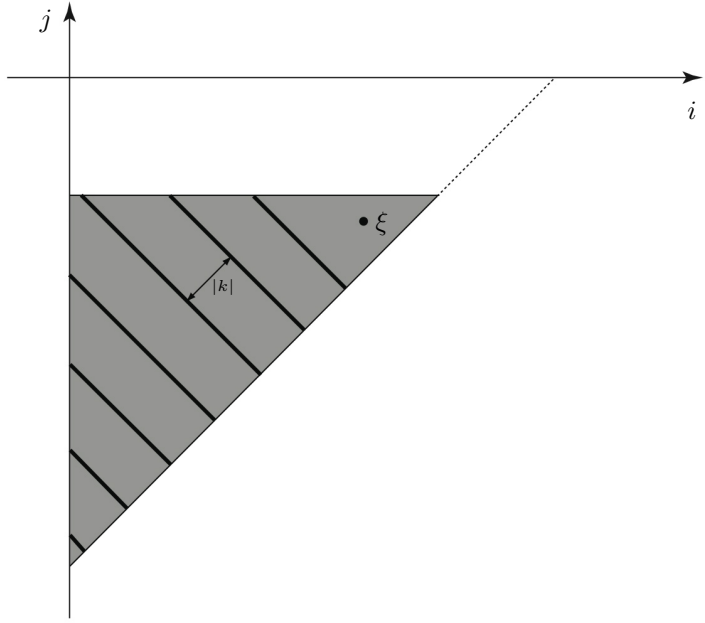}
\put(-15,215){$i$}
\put(-260,250){$j$}
\put(-133,173){$\xi$}
\put(-195,150){\tiny $|k|$}
\caption{The embedding $\varphi$. Letting the $ij$-plane represent $\Lambda ^*H^1(\Sigma _g)\otimes \mathbb{Z}[U,U^{-1}]$ as before, the shaded triangle represents $X(g,d)$ via the embedding \eqref{xgddef}. For an element $\xi \in X(g,d)$, represented by the black dot, $\varphi (\xi) - \xi$ is a sum of terms living on the lines
of slope $-1$ and of distance from $\xi$ increasing in multiples of $|k|$. Only lines to the \lq\lq left\rq\rq \,  of $\xi$ 
contain nontrivial summands of $\varphi (\xi) - \xi $. }  \label{pic5}
\end{figure}

\begin{theorem}\label{Znontorsresult} For any $k\neq 0$, there is an identification of vector spaces 
\[
HF^+(\Sigma\times S^1, \s_k) \cong X(g,d) \quad \quad (d = g-1-|k|).
\]
Under this identification, the action of $\Lambda^*H_1(\Sigma\times S^1)$ on $HF^+(\Sigma\times S^1,\s_k)$ is induced by the embedding $\varphi: X(g,d) \to \Lambda^*H^1(\Sigma)\otimes \zee[U^{-1}]$ given by 

\begin{equation}\label{phidef}
\varphi(\xi \otimes U^j) = \sum_{n\geq 0} (-1)^n \, (pr_{i\geq 0}\,U^{|k|} J^+)^n(\xi\otimes U^j)
\end{equation}
(see figure \ref{pic5}). When $3|k| > g-2$, this action can be identified with the standard one.
\end{theorem}
As remarked before, much of this theorem was obtained by Ozsv\'ath and Szab\'o in \cite{OSknot}: the group structure, and the action of $H_1(\Sigma\times S^1)$ when $3|k|>g-2$. We will see below that when the latter condition is not satisfied, the action is not standard in general.

\begin{proof} From the preceding lemma, to prove the first two statements it suffices to show that $\varphi$ induces an isomorphism $X(g,d)\to \kerr(F)$. Certainly $\varphi$ is injective, and the fact that $\imm(\varphi) \subset \kerr(F)$ is easy to check. On the other hand, by considering top-degree parts one can verify that any element of $\kerr(F)$ is of the form given on the right-hand side of \eqref{phidef}. Namely, let us pick  
$\xi \in \kerr F$ of the form $\xi = \xi_\ell+\xi_{\ell - 2|k|}+\xi_{\ell-4|k|}+ ...$ where $\xi_i$ is the homogeneous part of $\xi$ of degree $i$ and where we assume that $\xi_{\ell + 2n|k|}=0$ for all $n>0$.  
The equation $F(\xi) = (\mbox{id} + pr_{i\ge 0}U^{|k|}J^+)\xi = 0$ then becomes the system of equations 
\begin{align} \nonumber
\xi_ \ell & = 0 \cr
\xi_{\ell-2|k|} +pr_{i\ge 0} U^{|k|} J^+ (\xi _\ell) & = 0 \cr  
\xi_{\ell-4|k|} +pr_{i\ge 0} U^{|k|} J^+ (\xi _{\ell-2|k|}) & = 0 \cr  
& \, \, \, \vdots 
\end{align}
showing that $\xi_{\ell - 2n|k|} = (-1)^n (pr_{i\ge 0} U^{|k|}J^+) \xi _\ell $ for all $n\ge 0$ and so $\xi = \varphi (\xi _\ell)$. A general kernel element of $F$ is a sum of such $\xi$'s taken over congruency classes of $\ell \mod 2|k|$. It is easily seen that the above system of equations also splits into systems
according to the congruency class of $\ell \mod 2|k|$, thus proving that $\kerr F = \imm \varphi$. 

To prove that the action is standard for $3|k|>g-2$, consider the more general question of the determination of the action of $\gamma\in H_1(\Sigma)$ on a homogeneous element $\xi\in X(g,d)$, under the identification $\varphi$. That is to say, we must compute $\varphi^{-1}(\gamma.\varphi(\xi))$, where the action of $\gamma$ on $\Lambda^*H^1(\Sigma)\otimes \zee[U^{-1}]$ is given by \eqref{theaction}; in particular we must compare $\varphi^{-1}(\gamma.\varphi(\xi))$ with $\gamma. \xi$.

We can describe the action of $pr_{i\geq 0} U^{|k|}J^+$ diagrammatically as follows. An element at position $(i,j)$ is mapped into the group at position $(j,i)$, then terms are added in all groups on that diagonal lying below and to the right of this position (corresponding to the contraction with $\exp(2\omega U)$), but lying on or above the $i$ axis. Then all terms are pushed diagonally down and left by $|k|$ steps, and any portion to the left of the $j$ axis is killed by $pr_{i\geq 0}$. In particular, if $\xi\in C\{i,j\}$ is such that the result of this sequence of operations is zero (so the $n = 1$ and higher terms of \eqref{phidef} vanish), then the $H_1(\Sigma)$ action is standard on $\xi$. 

Now, observe that in $\varphi(\xi) = \xi - (pr_{i\geq 0}U^{|k|}J^+)(\xi) + (pr_{i\geq 0}U^{|k|}J^+)^2(\xi) - \cdots$, the only term lying in $C\{i\geq 0\mbox{ and } j\leq -|k|-1\}$ is the first, $\xi$ itself. Put another way, an element of $\kerr(F)$ is determined uniquely by its component in $\iota(X(g,d))$. From the diagrammatic description of $\varphi$ above, it is easy to see that the terms of the difference $\gamma.\varphi(\xi) - \gamma.\xi$ that lie in $\iota(X(g,d))$ are in the top row of the latter subgroup, namely $C\{*, -|k|-1\}$. Since elements of $C\{*,-|k|-1\}$ have degree at least $-|k|-1$, it follows that so long as
\[
\deg(\gamma.\xi) - 2|k| < -|k|-1,
\]
we will have $\iota^{-1}(\gamma.\xi) = \varphi^{-1}(\gamma.\varphi(\xi))$. The above is equivalent to $\deg(\xi) < |k|$, and since the maximum degree of an element of $X(g,d)$ is $g-2|k|-2$ the last claim of the theorem follows.
\end{proof}

The description given in the introduction of the corrections to the standard $H_1$ action appearing when $3|k|\leq g-2$ follows from the discussion in the preceding proof. Indeed, we saw above that the difference $\varphi^{-1}(\gamma.\varphi(\xi)) - \iota^{-1}(\gamma.\xi)$ is in general a combination of terms $\rho_\ell(\gamma,\xi)$ lying in $C\{*, -|k|-1\}$. The latter term is the component of $(pr_{i\geq 0}U^{|k|}J^+)^\ell(\gamma.\xi)$ that lies in $\iota(X(g,d)) = C\{i\geq 0\mbox{ and } j\leq -|k|-1\}$ (equivalently, the component of this quantity lying in the lowest available exterior power). This component is explicitly computable from the expression for $J^+$ in Proposition \ref{Jprop}.

\subsection{The case of $HF^+(Y_0, \s_0)$} We work with coefficients in $\cee$. The calculation of $HF^+$ in the case of $k = 0$, corresponding to the torsion \spinc structure on $\Sigma\times S^1$, follows much the same lines as the calculation of $\widehat{HF}$ in this case. As before we consider the surgery sequence for negative surgeries \eqref{LESplus}, identifying the map $F^+: HF^+(Y,\s)\to HF^+(Y_{-n},\t_0)$ with $v^+ + h^+$ where in this case we have
\[
h^+ = pr_{i\geq 0} \circ J^+.
\]

We claim that for any $\xi\in HFK^\infty$, we have
\begin{equation}\label{projection}
F^+(pr_{i\geq 0}\xi) = \pi_C(1 + J^\infty)(\xi),
\end{equation}
where $\pi_C = pr_{j\geq 0}\circ pr_{i\geq 0}$. To see this, observe that $F^+ = \pi_C \circ(1 + J^+)$, and since $J^\infty$ maps $\{i<0\}$ into $\{j<0\}\subset \ker(\pi_C)$, we have $J^+(pr_{i\geq 0} \xi) = pr_{j\geq0}J^\infty(\xi)$ for any $\xi\in HFK^\infty$.  

For the statement of the next lemma, recall that a form $\xi\in\Lambda^*H^1(\Sigma _g)$ is coprimitive if $\omega\wedge\xi = 0$.

\begin{lemma}\label{imagelemma} Let $V\subset \Lambda^*H^1(\Sigma)\otimes \zee[U,U^{-1}]$ denote the $\zee[U,U^{-1}]$ submodule given by
\[
V = \bigoplus_{\mbox{\scriptsize $\begin{array}{c}j = 0,\ldots,g-1 \\ k = 0,\ldots g-j-1 \end{array}$}} \omega^k P^j \otimes \zee[U,U^{-1}].
\]
Thus $V$ is complementary to the submodule spanned by the coprimitive forms. Then $\imm(1 + J^\infty) = e^{-\omega U}\wedge V$. 
\end{lemma}

\begin{proof} For notational convenience, we drop reference to $U$; since all maps involved preserve the absolute grading, the power of $U$ attached to a given homogeneous form in a particular absolute grading can be recovered from the grading and the degree of the form. 

Recall from Lemma \ref{altJformulalemma} that $J^\infty( e^{-\omega}\wedge\xi) = -e^{-\omega}\wedge (e^{-\omega}\angle\xi))$, so that if $\phi:\Lambda^*H^1(\Sigma)\to \Lambda^*H^1(\Sigma)$ is the isomorphism given by $\phi(\xi) = e^{-\omega}\wedge \xi$ then $(1+J)(\phi(\xi)) = \phi((1-e^{-\omega})\angle \xi)$. Therefore the lemma is equivalent to the fact that the image of contraction by $1-e^{-\omega}$ is $V$. To see this, recall that contraction by $\omega$ induces an isomorphism $\omega^k P^j\to \omega^{k-1} P^j$ for all $k\geq 1$. Hence, in particular, each $P^j$ with the exception of $P^g$ is in the image of contraction by $\omega$, and therefore contraction by $1-e^{-\omega}$ maps $\omega P^j$ onto $ P^j$ for all $j<g$. Continuing inductively on $k>1$, we see that given $x\in \omega^k P^j$, $k<g-j$, there exists $y\in \omega^{k+1} P^j$ with $\omega\angle y = x$. Therefore $(1-e^{-\omega})\angle y = x + $ lower order terms, where the lower order terms are in the image of contraction by $1-e^{-\omega}$ acting on $\omega^{<k}P^j$. In particular, $x$ is in $\imm((1-\exp(\omega))\angle\,\cdot\,)|_{\omega^{<k+1} P^j})$.
\end{proof}

We are interested in determining the kernel and cokernel of $F^+$, and since $F^+$ preserves degree, we may work in a given absolute grading $d$. For $d\geq g-1$ the restrictions of $pr_{i\geq 0}$ and $\pi_C$ to degree $d$ are the identity, so that $F^+ = 1+J^\infty$ in these degrees. In particular, if $d\geq g-1$, it follows from Corollary \ref{espacecor} and the lemma preceding it that 
\begin{equation} \label{kernelidentivo}
\kerr(F^+_d) = e^{-\omega U}\wedge ( P^*\cap HF^+_d(Y))
\end{equation}
where $ P^*\subset \Lambda^*H^1(\Sigma)\otimes \zee[U,U^{-1}]$ denotes the submodule spanned by the primitive forms in $\Lambda^*H^1(\Sigma)$, and $F^+_d$ is the restriction of $F^+$ to $HF^+_d(Y)$. Likewise, it follows from \eqref{projection} and Lemma \ref{imagelemma} that
\begin{equation}\label{cokerident}
\cokerr(F^+_d) \cong \tilde{ P}^*\cap HF^+_d(Y_{-n})
\end{equation}
where $\tilde{ P}^*\subset \Lambda^*H^1(\Sigma)\otimes \zee[U,U^{-1}]$ is the submodule spanned by coprimitives (note that $HF^+_d(Y_{-n})$ and $HF^+_d(Y)$ are isomorphic in this degree range). Note also that the coprimitive subspace of $\Lambda^*$ is just $\bigoplus_j \omega^{g-j}P^j$.

\begin{lemma}\label{cokerlemma} The cokernel of $F^+$ is given by \eqref{cokerident} for all $d\geq 0$.
\end{lemma}

\begin{proof}
We already know \eqref{cokerident} holds for $d\geq g-1$. Fix $d< g-1$ (we can assume $d\geq 0$ since $HF^+_d(Y_{-n})=0$ for $d<0$); then $F^+_d$ is identified with a map 
\[
F^+_d : \Lambda^{g+d} \oplus \Lambda^{g+d-2} \oplus\cdots\oplus \Lambda^{g-d}\oplus\Lambda^{g-d-2} \oplus\cdots \longrightarrow \Lambda^{g+d} \oplus \Lambda^{g+d-2} \oplus\cdots\oplus \Lambda^{g-d},
\]
namely the restriction of $\pi_C(1+J^\infty)$. The proof of Lemma \ref{imagelemma} shows that the image of $(1+J^\infty)|_{\Lambda^{\leq g+d}}$ is equal to $e^{-\omega}\wedge (\Lambda^{\leq g+d-2} \cap V_d)$, where $V_d$ is the subspace of $V$ lying in degree $d$ (in fact, it shows that $\imm((1+J^\infty)|_{e^{-\omega}\wedge \Lambda^{\leq r}}) = e^{-\omega}\wedge(\Lambda^{\leq r-2}\cap V_d)$, which is enough for us since $pr_{i\geq 0}(e^{-\omega}\wedge \Lambda^{\leq g+d}) = \Lambda^{\leq g+d}$ here). We will show that the projection of this space under $\pi_C$ is  $e^{-\omega} \wedge ((\Lambda^{g+d}\oplus\cdots\oplus \Lambda^{g-d}) \cap V_d)$:
\begin{equation}\label{toshow}
\pi_C(e^{-\omega}\wedge (V_d\cap\Lambda^{\leq g+d-2})) = e^{-\omega} \wedge ( V_d\cap (\Lambda^{g+d}\oplus\cdots\oplus \Lambda^{g-d})).
\end{equation}
Here and below it is understood that on the right hand side, multiplication by $e^{-\omega}$ occurs in $HF^+_d(Y_{-n}) = \Lambda^{g+d}\oplus\cdots\oplus\Lambda^{g-d}$, so that forms of degree greater than $g+d$ do not appear.

Define $S^j =  P^j \oplus \omega  P^j \oplus \omega^2  P^j \oplus\cdots\oplus \omega^{g-j-1} P^j$, so $S^j\cap \Lambda^{\leq g+d}$ is the part of $V\cap HF^+_d(Y)$ generated by primitive forms of degree $j$. That is, 
\[
 V_d\cap \Lambda^{\leq g+d} = \bigoplus_{j=0}^g S^j\cap \Lambda^{\leq g+d}
\]
To prove \eqref{toshow} it suffices to show
\begin{equation}\label{toshow2}
\pi_C(e^{-\omega}\wedge(S^j\cap \Lambda^{\leq g+d-2})) = e^{-\omega}\wedge(S^j\cap (\Lambda^{g+d}\oplus\cdots\oplus\Lambda^{g-d}))
\end{equation}
for each $j$. If $j\geq g-d$, we have that $S^j\subset \Lambda^{g-d}\oplus\cdots\oplus \Lambda^{g+d-2}$ and $e^{-\omega} S^j\subset \Lambda^{g-d}\oplus\cdots\oplus \Lambda^{g+d}$, so that $\pi_C|_{e^{-\omega} S^j} = \mbox{id}_{e^{-\omega} S^j}$. Hence \eqref{toshow2} holds for such $j$. 

Fix $j<g-d$, and suppose for convenience that $g-d-j=2k$ is even. If $\{b\}_{b\in \B}$ is a basis for $\omega^{k-1}P^j$, then a basis for $S^j\cap (\Lambda^{g-d} \oplus\cdots\oplus\Lambda^{g+d})$ is given by $\B' = \{\eta_1b, \eta_2b, \ldots, \eta_db, \eta_{d+1}b\}_{b\in \B}$. (Note that for this range of $j$, multiplication by $e^{-\omega}$ is an isomorphism of $S^j\cap(\Lambda^{g+d}\oplus\cdots\oplus\Lambda^{g-d})$ so may be ignored.) Hence the left hand side of \eqref{toshow2} contains the images under $\pi_C$ of $\{e^{-\omega} b, e^{-\omega}\eta_1b, e^{-\omega}\eta_2b, \ldots, e^{-\omega}\eta_db\}_{b\in\B}$. Since $e^{-\omega} = \sum_{j\geq 0}(-1)^j\eta_j$ and $\eta_j\eta_i = {{i+j}\choose{i}}\eta_{i+j}$, the expression of $\pi_C(e^{-\omega}\eta_ib)$ in terms of the basis $\B'$ is
\[\textstyle
\pi_C(e^\omega\eta_i b) = \left((-1)^{i+1}{{1}\choose{i}}\eta_1b, (-1)^{i+2}{{2}\choose{i}}\eta_2b,\ldots, (-1)^{i+d+1}{{d+1}\choose{i}}\eta_{d+1}b\right),
\]
where ${{n}\choose{k}} = 0$ if $n<k$. Letting $i$ vary from $0$ to $d$, we have that the left hand side of \eqref{toshow2} contains ($\dim( P^j)$ copies of) the column space of the matrix
\[
\left[\begin{array}{cccccccc} -1 & 1 \\ 1 & -2 & 1 \\ -{{3}\choose{0}} & {{3}\choose{1}}& -{ 3\choose 2} & 1 \\ \\& \vdots & &&& \ddots \\\\ (-1)^d{{d}\choose 0} & (-1)^{d+1}{{d}\choose 1} & (-1)^{d+2}{d\choose 2} & & \cdots & & -{d\choose{d-1}} & 1 \vspace{.5ex}\\ (-1)^{d+1}{{d+1}\choose 0} & (-1)^{d+2} {{d+1}\choose 1} & (-1)^{d+3}{{d+1}\choose 2} && \cdots && {d+1}\choose{d-1} & -{{d+1}\choose d} \end{array}\right]
\]
Successively adding the first column to the second, the second to the third, etc., it is easy to see that the determinant of this matrix is $(-1)^{d+1}$, hence it is surjective. This proves \eqref{toshow2} in case $g-d-j$ is even; the odd case is similar, and \eqref{toshow} is proved. 

Finally, we observe that since $V$ is complementary to the coprimitive forms, and multiplication by $e^{-\omega}$ is an isomorphism of $HF^+_d(Y_{-n})$ preserving the coprimitives, the fact that the image of $F^+_d$ is given by the right hand side of \eqref{toshow} proves \eqref{cokerident}.
\end{proof}

As we shall see in the proof of theorem \ref{noviteorem}, the analogue of Lemma \ref{cokerlemma} does not hold for $\kerr(F^+)$ (i.e. equation \eqref{kernelidentivo} does not hold for arbitrary $d$) . 
However, from the present results, we do obtain a concrete description of $HF^+(\Sigma_g\times S^1,\s_0; \cee)$ in large degrees---equivalently, a description of $HF^\infty(\Sigma_g\times S^1,\s_0;\cee)$. In general, 
\[
HF^+_{d+1/2}(\Sigma_g\times S^1) \cong \kerr(F^+_d) \oplus \cokerr(F^+_{d+1}).
\]
In a given degree $d\geq g-1$, $\kerr(F^+_d)$ is isomorphic to the primitive forms in even (or odd) degree while $\cokerr(F^+_{d+1})$ is isomorphic to the (co)primitive forms in odd (or even) degree. Hence we have an isomorphism 
\[
HF^+_{d+1/2}(\Sigma_g\times S^1) \cong  P^* \qquad (d\geq g-1)
\]
of the Floer homology of $\Sigma_g\times S^1$ with the primitive cohomology of the Jacobian torus $T^{2g}$, for all $d\geq g-1$. In particular, we see that $HF^\infty(\Sigma_g\times S^1,\s_0;\cee)$ has dimension ${{2g+1}\choose{g}}$ in each degree.

\begin{theorem} \label{noviteorem}
There is an isomorphism of graded $\zee[U]$-modules
\begin{eqnarray}
HF^+(\Sigma_g\times S^1, \s_0; \cee) &\cong& HF^+_{red}(\Sigma_g\times S^1, \s_0;\cee) \oplus \nonumber\\&& \hspace*{.5in}\bigoplus_{j\geq 0} (P^j\otimes \T_{-g+j+1/2}) \oplus \bigoplus_{j\geq g} \tilde{P}^j\otimes \T_{-g+j - 1/2}.\label{theanswer}
\end{eqnarray}
Here $\T_n$ denotes the module $\zee[U,U^{-1}]/U\cdot \zee[U]$, graded so that the smallest degree of a nonzero homogeneous element is $n$, while $\tilde{P}^j$ denotes the space of coprimitive forms in $\Lambda^jH^1(\Sigma_g)$.

Furthermore, there is an isomorphism of graded abelian groups
\begin{equation}\label{redanswer}
HF^+_{red}(\Sigma_g\times S^1,\s_0;\zee) \cong X_0(g,g-3)[{\ts\frac{5}{2}}],
\end{equation}
where $X_0(g,d)$ denotes the group underlying the $\zee[U]$-module $X(g,d)$ (we take $X(g,d) = 0$ for $d<0$).
\end{theorem}

Thus there is no reduced Floer homology for $g = 1, 2$, and when $g\geq 3$ there is nontrivial reduced Floer homology in degrees between $-g+\frac{5}{2}$ and $g-\frac{7}{2}$. We remark that the identification of the reduced Floer homology with $X(g,g-3)$ does not respect the $U$-action.

\begin{proof} If $x = pr_{i\geq 0}(e^{-\omega U}\beta\otimes U^q)\in HF^+_d(Y)$, for $\beta$ a primitive form and $d<g-1$, then $x$ is the image under $U^n$ of $e^{-\omega U}\beta \otimes U^{q-n}\in HF^+_{d+2n}(Y)$. Likewise, it follows from \eqref{projection} and Corollary \ref{espacecor} that all elements of $\kerr(F^+)$ lying in the image of a large power of $U$ are of this form. Similarly, since the cokernel of $F^+$ is given by (the projection under $\pi_C$ of) the coprimitive forms, all cokernel elements are ``reducible'' in this sense. Thus for a given degree $d<g-1$, Floer homology classes in $HF^+_{d+1/2}(\Sigma\times S^1,\s_0;\cee)$ lying in the image of a large power of $U$ are identified with $\pi_{i\geq 0}(e^{-\omega U}\wedge\beta)$ for a primitive form $\beta$, or $\pi_C(\tilde{\beta})$ for $\tilde{\beta}$ coprimitive. Observe that the kernel of the action of $U$ on $HF^+_d(Y)$ is $\Lambda^{g+d} = C\{0,d\}$, while on $HF^+_d(Y_{-n})$ it is $\Lambda^{g+d}\oplus \Lambda^{g-d}$ (when $d>0$). Since there are primitive forms only in exterior powers of degree $\leq g$, ($e^{\omega U}$ times) primitive forms are killed by $U$ once $d\leq 0$: specifically, $e^{-\omega U} P^j$ is mapped to $0$ by $U$ acting in degree $-g+j$. This accounts for the factor $ P^j\otimes \T_{-g+j+1/2}$ in \eqref{theanswer}. Similarly, coprimitive forms lie in exterior powers of degree $\geq g$, and $\tilde{ P}^j$ is killed by the action of $U$ (on $HF^+(Y_{-n})$) in degree $-g+j$ with $j\geq g$. The connecting homomorphism in the surgery sequence decreases degree by $\frac{1}{2}$, so we see the factor of $\tilde{ P}^j\otimes \T_{-g+j-1/2}$ in \eqref{theanswer}. Furthermore, we have accounted for all non-reduced Floer homology as well as all of the cokernel of $F^+$ by Lemma \ref{cokerlemma}. It remains to identify the portion of the kernel complementary to that given by primitive forms with $X_0(g,g-3)[\frac{5}{2}]$.

We do this first on the level of graded vector spaces (i.e., with coefficients in $\cee$) by a simple dimension count.

For negative degrees, the calculation is particularly easy. Identify the triangle $C\{i\geq 0\mbox{ and } j<0\}$ with $X(g,g-1)$ as previously, and note that the kernel of $F^+$ in negative degrees is the portion of this triangle lying below the line $i=-j$. Therefore, we are interested in the ``non-primitive'' portion of $X(g,g-1)$ (note that ``primitive'' here does not refer to the Lefschetz decomposition of $H^*(\mbox{Sym}^{g-1}\Sigma_g)$). For $n\leq g$, we have $\Lambda^n/ P^n \cong\Lambda^{n-2}$, so that $X(g,g-1)$ modulo primitives is isomorphic to $X(g,g-3)[2]$ (c.f. \eqref{xgddef}). Taking account of the degree shift in the surgery sequence gives \eqref{theanswer} for negative degrees.

In general, we have the exact sequence
\[
0\to \cokerr(F^+_{d+1})\to HF^+_{d+\frac{1}{2}}(\Sigma_g\times S^1)\to HF^+_d(Y)\to HF^+_d(Y_{-n})\to \cokerr(F^+_{d})\to 0,
\]
where we can assume $0\leq d < g-1$. Write $c_d = \dim(\cokerr(F^+_d))$, $a_d = \dim(HF^+_{d}(Y))$, and $b_d = \dim(HF^+_{d}(Y_{-n}))$. Then
\[
\dim(HF_{d+\frac{1}{2}}^+(\Sigma_g\times S^1)) = c_{d+1} + c_{d} + a_d - b_d
\]
while
\[
\dim(\kerr(F^+_d)) = a_d - b_d + c_d.
\]
The ``reducible part'' of the kernel (i.e., the part in the image of a large power of $U$) is identified with the primitives lying in $HF^+_d(Y)$. Since $d\geq 0$ this includes all the primitives whose degree is of the appropriate parity, thus
\[
\dim(\mbox{reducibles}) = \left\{\begin{array}{ll} {{2g}\choose{g}} & g-d\equiv 0\mod 2\vspace{1ex}\\ {{2g}\choose{g-1}} & g-d\equiv 1 \mod 2\end{array}\right.
\]
The dimension of $HF^+_{red, d+\frac{1}{2}}(\Sigma_g\times S^1)$ is given by the difference of the preceding two expressions. Now, it is easy to see
\[
a_d - b_d = {{2g}\choose{g-d-2}} + {{2g}\choose{g-d-4}} + \cdots,
\]
while from Lemma \ref{cokerlemma} we have
\begin{eqnarray*}
c_d &=& \dim(\tilde{P}^{g+d} \oplus \tilde{P}^{g+d-2} \oplus \cdots) \\&=& \ts\left({{2g}\choose{g-d}} - {{2g}\choose{g-d-2}}\right) + \left({{2g}\choose{g-d+2}} - {{2g}\choose{g-d}}\right) + \cdots \\&=& \left\{\begin{array}{ll} {{2g}\choose{g}} - {{2g}\choose{g-d-2}} & g-d\equiv 0\mod 2\vspace{1ex} \\ {{2g}\choose{g-1}} - {{2g}\choose{g-d-2}} & g-d\equiv 1 \mod 2\end{array}\right.
\end{eqnarray*}
Hence 
\begin{eqnarray*}
\dim HF^+_{red,\,d+\frac{1}{2}}(\Sigma_g\times S^1) &=& a_d-b_d+c_d - \dim(\mbox{reducibles})\\
&=& \ts{{2g}\choose{g-d-4}} + {{2g}\choose{g-d-6}} + \cdots\\
&=& \dim(X_{d-2}(g,g-3)),
\end{eqnarray*}
which confirms the structure of $HF^+_{red}(\Sigma_g\times S^1)$ as a graded vector space.

We now need only see that the reduced Floer homology is free over $\zee$ to verify \eqref{redanswer}. This is a simple consequence of the facts that $HF^+_{red}(Y_{-n};\zee)=0$, hence $HF^+_{red}(\Sigma\times S^1; \zee)\to HF^+(Y;\zee)$ is injective; and that $HF^+(Y;\zee)$ is free. 
\end{proof}

\begin{remark}\label{remarks} The action of $U$ on $HF^+_{red}(\Sigma\times S^1;\s_0)$ can be read off in part from the long exact sequence \eqref{Useq} relating $HF^+$ and $\hfhat$. Indeed, it is a straightforward exercise to verify the following:
\begin{enumerate}
\item $U:HF^+_{red, d}(\Sigma\times S^1)\to HF^+_{red, d-2}(\Sigma\times S^1)$ is surjective whenever $d$ is less than or equal to the middle nontrivial grading; i.e., for $d\leq -\frac{1}{2}$. 
\item $U:HF^+_{red, d}(\Sigma\times S^1)\to HF^+_{red, d-2}(\Sigma\times S^1)$ is injective above the middle degree, i.e., for $d>\frac{3}{2}$. 
\item $U: HF^+_{red, \frac{3}{2}}(\Sigma\times S^1)\to HF^+_{red, -\frac{1}{2}}(\Sigma\times S^1)$ has kernel of dimension 
\[\ts
2^{g-1} - \frac{1}{2}{{2g}\choose{g}} + {{2g}\choose{g-2}}.
\]

\end{enumerate}
One can identify the ``unexpected'' kernel elements of $U$ in (3) above with (the preimage in $HF^+(\Sigma\times S^1)$ of) the non-primitive $\tstar$-self-dual vectors lying in $C\{0,0\} = \Lambda^g$ (c.f. section \ref{hatsubsec}). Note that the action of $U$ on $X(g,g-3)$ in the corresponding degree is an isomorphism, in contrast to the situation here.

Furthermore, since $HF^+_d(\Sigma\times S^1)\to HF^+_{d-\frac{1}{2}}(Y)$ is an isomorphism for $d<-\frac{1}{2}$, we can say that $HF^+_{red}(\Sigma\times S^1,\s_0;\zee)$ is identified with $X(g,g-3)$ as a $U$-module in those degrees. 
\end{remark}

\begin{remark} We have seen that $HF^+_{d+1/2}(\Sigma\times S^1,\s_0;\zee)$ is isomorphic to $\kerr(F^+_d)\oplus \cokerr(F^+_{d+1})$, and in large degrees we can identify $F^+_d$ with $1 + J^\infty$. Furthermore, as in the proof of Lemma \ref{imagelemma}, after multiplication by $e^{-\omega U}$ we can identify $1+J^\infty$ with contraction by $1-e^{-\omega U}$. The kernel of this homomorphism is given by the primitive forms; hence the determination of $HF^\infty(\Sigma\times S^1,\s_0;\zee)$ reduces to the calculation of the cokernel of contraction by $1-e^{-\omega}$ acting on $\Lambda^*H^1(\Sigma)$. The latter is an explicit but complicated calculation that can be carried out by computer for small values of $g$.

It is interesting to note that the cokernel of contraction by $\omega$ itself has been computed in general by Lee and Packer \cite{leepacker} in another context, and in all cases so far computed yields the same result as contraction with $1-e^{-\omega}$. In particular, it can be seen that for all $g\leq 7$, we have $HF^\infty(\Sigma_g\times S^1,\s_0;\zee) \cong H^*(E_g;\zee)\otimes \zee[U,U^{-1}]$ after a shift in grading, where $E_g$ is the $S^1$-bundle over the Jacobian torus of $\Sigma_g$ discussed in the introduction.
\label{contrrem}
\end{remark}

\subsection{Calculation modulo 2}\label{mod2sec} The above calculation of Floer homology for $\Sigma\times S^1$ in the torsion \spinc structure is valid only for coefficients in the complex numbers, mainly because of the central use of the decomposition of $\Lambda^*H^1(\Sigma)$ into $\SL_2(\cee)$ representations. We cannot use this decomposition in general, but in the case of coefficients in $\zee/2\zee$ the calculation proceeds fairly easily. 

Observe first that the two ``known'' terms in the surgery sequence \eqref{LESplus} are free abelian groups, so we pass to $\zee_2$ coefficients merely by forming the tensor product of all groups and maps with $\zee_2$. Furthermore, a glance at Proposition \ref{Jprop} shows that when working with coefficients mod 2, our isomorphism $J^+$ is identified with the Hodge-Lefschetz operator $\tstar$. Thus, we are interested in the kernel and cokernel of
\begin{eqnarray*}
F^+: C\{i\geq 0\} & \to & C\{i\geq 0 \mbox{ and } j\geq 0\}\\
\xi\otimes U^n &\mapsto & \pi_C(\xi\otimes U^n + \tstar \xi\otimes U^m)
\end{eqnarray*}
where $\pi_C$ is the projection to $C\{i\geq 0 \mbox{ and } j\geq 0\}$ as before, and $m = g+n-p$ if $\xi\in \Lambda^pH^1(\Sigma)$. 

From this description, we can understand $F^+$ by considering its restrictions to the summands $C\{i,i\}$ and $C\{i,j\}\oplus C\{j,i\}$ for $i\neq j$. It is clear that if $\xi \in C\{i,j\}$ with $j< 0$ then $F^+(\xi) = 0$. If $i,j\geq 0$ and $i\neq j$, then the restriction of $F^+$ to $C\{i,j\}\oplus C\{j,i\}$ appears modulo signs as
\[
F^+ = \left[\begin{array}{cc} 1 & \tstar \\ \tstar & 1\end{array}\right]
\]
Since $\tstar$ is an isomorphism whose square is (plus or minus) the identity, this means that the kernel and cokernel of the above map are each isomorphic to $C\{i,j\}\cong \Lambda^pH^1(\Sigma)$---in particular, each has dimension ${2g}\choose{p}$. Furthermore, the action of powers of $U$ maps elements of the form $(x, \tstar x)\in C\{i,j\}\oplus C\{j,i\}$ onto any desired element $U^n x$ of $C\{i\geq 0\mbox{ and } j<0\}$. Likewise, all representatives of the cokernel of $F^+$ are in the image of large powers of $U$.

Finally we consider the action of $F^+$ on $C\{i,i\}\cong \Lambda^gH^1(\Sigma)$. The calculation here proceeds just as in the calculation of $\hfhat(\Sigma\times S^1,\s_0;\zee)$ in section \ref{hatsubsec} and shows that the kernel of $F^+ = 1 + \tstar$ is generated by elements of the form 
\[
x_{i_1}\cdots x_{i_p}(z_{j_1}\cdots z_{j_q} + z_{k_1}\cdots z_{k_q}).
\]
Clearly such elements are carried to one another by powers of $U$. Thus the reduced Floer homology with coefficients in $\zee_2$ is trivial, and comparing with section \ref{hatsubsec} we obtain:

\begin{theorem} There is an isomorphism of graded $\zee_2[U]$-modules
\[
HF^+(\Sigma\times S^1,\s_0; \zee_2) \cong \hfhat(\Sigma\times S^1,\s_0;\zee)\otimes \T_0\otimes\zee_2.
\]
\end{theorem}

In particular, we see that the rank of $HF^\infty(\Sigma\times S^1,\s_0;\zee_2)$ in any fixed degree is equal to half the total rank of $\hfhat(\Sigma\times S^1,\s_0;\zee)$, namely $2^{2g-1} + 2^{g-1}$, while according to previous work that of $HF^\infty(\Sigma\times S^1,\s_0;\cee)$ is ${2g+1}\choose{g}$. Since the former is strictly greater than the latter whenever $g > 2$, we have:

\begin{corollary} For all genera $g >2$, the Heegaard Floer homology $HF^\infty(\Sigma_g\times S^1, \s_0; \zee)$ contains non-trivial 2-torsion.
\end{corollary}


\begin{thebibliography}{99}

\bibitem{GH} Phillip Griffiths and Joseph Harris, Principles of Algebraic Geometry, John Wiley \& Sons, 1978.
\bibitem{us1} Stanislav Jabuka and Thomas Mark, Heegaard Floer homology of certain mapping tori, Algebr. Geom. Topol. 4 (2004) 685--719.
\bibitem{us2} Stanislav Jabuka and Thomas Mark, Heegaard Floer homology of mapping tori II in: H. Boden et. al. (Eds.), Geometry and topology of manifolds, Fields Inst. Commun. 47 Amer. Math. Soc., Providence, RI, 2005, pp. 119--135
\bibitem{us3} Stanislav Jabuka and Thomas Mark, Product formulae for Ozsv\'ath-Szab\'o 4-manifold invariants, arXiv:0706.0339.
\bibitem{km} Peter Kronheimer and Tomasz Mrowka, Monopoles and three-manifolds. In preparation. 
\bibitem{leepacker} Soo Teck Lee and Judith A. Packer, The cohomology of the integer Heisenberg groups, J. Algebra 184 (1996), no. 1, 230--250.
\bibitem{muwa} Vincente Mu\~noz and Bai-Ling Wang, Seiberg-Witten-Floer homology of a surface times a circle for non-torsion spin$^\mathbb{C}$ structures, Math. Nachr., 278 (2005), no. 7-8, 844--863.
\bibitem{OS1} Peter Ozsv\'ath and Zolt\'an Szab\'o, Holomorphic disks and topological invariants for closed three-manifolds, Ann. of Math. 159 (2004), no. 3, 1027--1158.
\bibitem{OS2} Peter Ozsv\'ath and Zolt\'an Szab\'o, Holomorphic disks and three-manifold invariants: properties and applications, Ann. of Math. 159 (2004), no. 3, 1159--1245.
\bibitem{OS3} Peter Ozsv\'ath and Zolt\'an Szab\'o, Holomorphic triangles and invariants for smooth four-manifolds, Adv. Math. 202 (2006), no. 2, 326--400.
\bibitem{OS4} Peter Ozsv\'ath and Zolt\'an Szab\'o, Absolutely graded Floer homologies and intersection forms for four-manifolds with boundary, Adv. Math. 173 (2003), no. 2, 179--261.
\bibitem{OSsymp} Peter Ozsv\'ath and Zolt\'an Szab\'o,
Holomorphic triangle invariants and the topology of symplectic
four-manifolds, Duke Math. J. 121 (2004), no. 1, 1--34.
\bibitem{OSknot} Peter Ozsv\'ath and Zolt\'an Szab\'o, Holomorphic disks and knot invariants, Adv. Math. 186 (2004), no. 1, 58--116.
\bibitem{OSplumb} Peter Ozsv\'ath and Zolt\'an Szab\'o, On the Floer homology of plumbed three-manifolds, Geom. Topol. 7 (2003), 185--224.
\bibitem{OSsurg} Peter Ozsv\'ath and Zolt\'an Szab\'o, Knot Floer homology and integer surgeries, preprint arxiv:math.GT/0410300.
\bibitem{rasmussen} Jacob Rasmussen, Floer homology of surgeries on two-bridge knots, Algebr. Geom. Topol. 2 (2002), 757--789.
\bibitem{roberts} Lawrence P. Roberts, Rational blow downs in Heegaard-Floer homology, preprint arxiv:math.GT/0607675.
\bibitem{salamon} Dietmar Salamon, Seiberg-Witten invariants of mapping tori, symplectic fixed points, and Lefschetz numbers, in: Proceedings of 6th G\"okova Geometry-Topology Conference, Turkish J. Math. 23 (1999), no. 1, 117--143.



\end{thebibliography}
\end{document}